\documentclass[reqno, a4paper, 11pt]{amsart}
\usepackage{amssymb, amsthm, amsmath}
\numberwithin{equation}{section}
\renewcommand{\d}{\mathrm{d}}
\newcommand{\x}{\ma{x}}
\newcommand{\w}{\ma{w}}
\newcommand{\y}{\ma{y}}
\newcommand{\z}{\ma{z}}
\renewcommand{\v}{\ma{v}}
\renewcommand{\u}{\ma{u}}
\renewcommand{\t}{\ma{t}}

\newcommand{\M}{\mathfrak{M}}
\renewcommand{\ss}{\mathfrak{S}}
\newcommand{\m}{\mathfrak{m}}

\newcommand{\mla}{\boldsymbol{\lambda}}
\newcommand{\mbeta}{\boldsymbol{\beta}}

\newcommand{\prop}{Proposition }

\newcommand{\TT}{\mathcal{T}}
\newcommand{\SSS}{\mathcal{S}}
\newcommand{\WW}{\mathcal{W}}
\newcommand{\R}{\mathbb{R}}
\newcommand{\F}{\mathbb{F}}

\newcommand{\Z}{\mathbb{Z}}
\newcommand{\N}{\mathbb{N}}
\newcommand{\Q}{\mathbb{Q}}
\newcommand{\cQ}{\overline{\mathbb{Q}}}

\newcommand{\bfP}{\mathbb{P}}
\newcommand{\A}{\mathbb{A}}

\newcommand{\ma}{\mathbf}

\newcommand{\mcal}{\mathcal}

\renewcommand{\le}{\leqslant}
\renewcommand{\ge}{\geqslant}
\renewcommand{\leq}{\leqslant}
\renewcommand{\geq}{\geqslant}

\newcommand{\ve}{\varepsilon}

\newcommand{\al}{\alpha}
\newcommand{\D}{\Delta}
\newcommand{\del}{\delta}

\newcommand{\be}{\beta}
\newcommand{\la}{\lambda}
\newcommand{\sfl}{\mathsf{\Lambda}}

\newcommand{\lab}{\label}

\newcommand{\h}{\mathbf{h}}

\newtheorem{thm}{Theorem}
\newtheorem*{thm*}{Theorem}
\newtheorem{lem}{Lemma}
\newtheorem{pro}{Proposition}

\newtheorem*{cor*}{Corollary}

\DeclareMathOperator{\rank}{rank}
\DeclareMathOperator{\proj}{Proj}
\DeclareMathOperator{\spec}{Spec}

\DeclareMathOperator{\sing}{sing}
\DeclareMathOperator{\meas}{meas}

\theoremstyle{definition}
\newtheorem*{ack}{Acknowledgement}
\newtheorem*{notat}{Notation}

\newcommand{\barray}{\begin{eqnarray*}}
\newcommand{\earray}{\end{eqnarray*}}
\newcommand{\hcf}{\mathrm{gcd}}
\renewcommand{\mod}{\hspace{-1mm}\pmod}
\newcommand{\colt}[2]{\genfrac{}{}{0pt}{1}{#1}{#2}}
\newcommand{\tstack}[3]{\colt{#1}{\colt{#2}{#3}}}

\newcommand{\ep}{\varepsilon}

\newcommand{\Xns}{X_{\mathrm{ns}}}

\begin{document}
\title{Rational points on quartic hypersurfaces}

\author{T.D. Browning}
\author{D.R. Heath-Brown}
\address{School of Mathematics, University of Bristol, Bristol BS8 1TW}
\email{t.d.browning@bristol.ac.uk}
\address{Mathematical Institute,
24--29 St. Giles', Oxford OX1 3LB}
\email{rhb@maths.ox.ac.uk}

\date{\today}

\begin{abstract}
Let $X$ be a projective  
non-singular quartic hypersurface of dimension $39$ or more,
which is defined over $\Q$. We show that $X(\Q)$ is non-empty provided
that $X(\R)$ is non-empty and $X$ has $p$-adic points for every prime $p$.
\end{abstract}

\maketitle

\tableofcontents

\section{Introduction}

Let $X\subset\bfP_\Q^{n-1}$ be a geometrically integral quartic
hypersurface defined over $\Q$. Thus we may suppose that $X$ is defined by an 
absolutely irreducible quartic form $F\in\Z[x_1,\ldots,x_n]$.
We will write $\sing(X)$ for the singular locus of $X$,
a projective subvariety of $X$, whose dimension is an integer in the
interval $[-1,n-3]$. 
The primary aim of this paper is to establish conditions on $X$ under 
which we can ensure that the set $X(\Q)$ is non-empty. Specifically, 
we would like to establish
the Hasse principle for a large class of quartic forms defined over $\Q$.
This states that in order for $X(\Q)$ to be non-empty it is necessary
and sufficient
that $X(\A_\Q)=X(\R)\times \prod_{p}X(\Q_p)$ is non-empty, where 
$\A_\Q$ is the set of ad\`eles on $\Q$.

That one cannot hope for the Hasse principle to hold for all quartic
hypersurfaces is demonstrated by the example
$$
X_1: \quad 4x_1^4+9x_2^4=8(x_3^4+x_4^4).
$$
It has been shown by Swinnerton-Dyer \cite{swd} that $X_1(\Q)$ is
empty, despite the fact that $X_1(\A_\Q)$ is non-empty.
This example is explained by the Brauer--Manin obstruction, the
Brauer set $X_1(\A_\Q)^{\mathrm{Br}}$ being empty in this 
instance. Colliot-Th\'el\`ene \cite[Appendix]{ct} has shown that 
the Brauer--Manin obstruction is void for non-singular quartic
hypersurfaces $X\subset \bfP_\Q^{n-1}$, with $n\geq 5$. Thus it is
natural to ask whether the Hasse principle holds for all such
hypersurfaces. When the underlying form takes an appropriate shape we can get
reasonably close to this prediction. Suppose, for example, that $X$ is
diagonal and non-degenerate. Then a classical application of the 
Hardy--Littlewood circle method (see \cite[Chapter 8]{dav-book}) 
will establish the Hasse principle for $n\geq 17$.
In fact, the $p$-adic conditions hold automatically for $n\geq 17$
when $X$ is diagonal, as shown by Davenport and Lewis \cite{DL}. Thus
all one needs to check is that the coefficients are not all of the
same sign.  The problem of establishing the Hasse principle for
general non-singular quartic forms is substantially harder.

Let us write $\Xns=X\setminus \sing(X)$ for the locus of
non-singular points on $X$.  Then there is a rather long-standing
result due to Birch \cite{birch} which establishes the existence of
$\Q$-rational points on $X$, under the assumption that 
$\Xns(\A_\Q)$ is non-empty and 
\begin{equation}\label{eq:50}
n-\dim \sing(X) \geq 50.
\end{equation}
In particular, this confirms the Hasse principle for hypersurfaces
defined by non-singular quartic forms over $\Q$ in at least $49$
variables. Birch goes even further and provides an asymptotic formula
for the number $N_{X}(P)$ of rational points $x\in X(\Q)$ whose 
height $H(x)$ is
bounded by $P$, as $P\rightarrow \infty$.  
Under the assumption that \eqref{eq:50} holds, this estimate takes the shape
\begin{equation}
  \label{eq:manin}
N_X(P)=c_X P^{n-4}\big(1+o(1)\big),
\end{equation}
and confirms the conjecture of Manin \cite{man} for this
particular family of hypersurfaces.  The constant
$c_X\geq 0$
is a product of local densities whose positivity can be established
under the assumption that $\Xns(\A_\Q)$ is non-empty.

Birch's seminal work has since been revisited and generalised in a 
number of different
ways.  Define $h=h(X)$ to be the least positive integer such that the 
quartic form
$F$ can be written identically as 
$$
A_1B_1 +\cdots+A_hB_h,
$$
for forms $A_i,B_i\in\Z[x_1,\ldots,x_n]$ of positive degree.
Then Schmidt \cite{schmidt} has shown that the asymptotic formula 
\eqref{eq:manin}
holds when $h(X)\geq 18432$.
In a rather different direction, Birch's main result has
been generalised to arbitrary number fields by
Skinner \cite{skinner}. In fact, Skinner also establishes weak approximation
for non-singular quartic hypersurfaces, when $n\geq 49$. 
In the present paper our goal is to
extend the admissible range of $n$ for which the Hasse principle 
holds, as follows.

\begin{thm}\label{main}
Let $X\subset \bfP_{\Q}^{n-1}$ be a quartic hypersurface, with
$$
n-\dim \sing(X) \geq 42.
$$
Assume that $\Xns(\A_\Q)$ is non-empty.
Then there exist constants $P_0\geq 1$ and $c>0$, such that 
$
N_X(P) \geq c P^{n-4}
$
for $P\geq P_0$. 
\end{thm}

In view of \eqref{eq:50}, we have therefore been able to save $8$
variables over the approach taken by Birch. 
Suppose that $X\subset \bfP_{\Q}^{n-1}$ is defined by a non-singular quartic
form in $n\geq 41$ variables. Then it follows from Theorem \ref{main}
that the Hasse principle holds for $X$. 
It seems very likely that a suitable modification of the argument would 
yield weak approximation, and it would be interesting to see whether our 
 main result could be generalised to the number field setting. 
With more work, it should be possible to replace the lower bound for
$N_X(P)$ with an appropriate asymptotic formula.
As in the work of Birch, we will use the Hardy--Littlewood circle
method to establish Theorem \ref{main}.  
We shall give an overview of the proof in \S \ref{section:over}.

One significant difference between our work and the work of Birch is 
in the treatment of the singular series 
\begin{equation}
  \label{eq:ss}
\ss:=\sum_{q=1}^\infty \frac{1}{q^n}\sum_{\colt{a=1}{\hcf(a,q)=1}}^q
S_{a,q},
\end{equation}
where  
\begin{equation}\label{eq:aq}
S_{a,q}:=\sum_{\x\bmod q}e^{2\pi i aF(\x)/q},
\end{equation}
which may or may not converge. 
For the case $d=4$, Birch only established the absolute convergence of 
$\ss$ under the assumption that \eqref{eq:50} holds.  
We are able to do rather better than this and will establish the following 
result in \S \ref{s:ss}. 

\begin{thm}\label{main-ss}
Let $X\subset \bfP_{\Q}^{n-1}$ be a quartic hypersurface, with
$$
n-\dim \sing(X) \geq 27.
$$
Then the singular series $\ss$ is absolutely convergent.
\end{thm}

It seems possible that one has absolute convergence as soon as 
$$
n-\dim \sing(X)\geq 6,
$$ 
but we are clearly a long way from proving this. 
Birch's theorem applies more generally to arbitrary hypersurfaces
$V\subset \bfP_\Q^{n-1}$ of degree $d\geq 3$. The outcome of his
investigation is that the clean Hasse principle holds when 
$$
n-\dim \sing(V)\geq 2+(d-1)2^d,  
$$
by which we mean that $V(\Q)$ is non-empty as soon as
$V_{\mathrm{ns}}(\A_\Q)$ is non-empty. 
It seems possible that the ideas contained in the present
paper could be adapted to obtain this same conclusion for a different
range of $n-\dim \sing(V)$. 
However, preliminary investigations suggest that aside from additional
difficulties intrinsic in handling forms of higher degree, this will only allow
us to replace $(d-1)2^d$ by $d(d+1)2^{d-3}$. This is patently weaker
for $d\geq 6$.

\begin{notat}
Throughout our work $\N$ will denote the set of positive
integers.  For any $\al\in \R$, we will follow common convention and
write $e(\al):=e^{2\pi i\al}$ and $e_q(\al):=e^{2\pi i\al/q}$. The
parameter $\ve$ will always denote a small positive real
number, which is allowed to take different values at different parts
of the argument.  We shall use $|\x|$ to denote the norm $\max |x_i|$ 
of a vector $\x=(x_1,\ldots,x_n)$. 
All of the implied constants that appear in this 
work will be allowed to
depend upon the coefficients of the quartic form $F$ under
consideration, the number $n$ of variables  involved, and the
parameter $\ve>0$.  Any further dependence will be explicitly
indicated by appropriate subscripts.
\end{notat}

\begin{ack}
This work was begun while the authors participated in
the programme ``Rational and integral points on higher-dimensional
varieties'' held at M.S.R.I.,
during the period 09/01/06-- 19/05/06.
The hospitality and financial support of the institute is gratefully 
acknowledged.
\end{ack}

\section{Overview of the proof}\label{section:over}

Our proof of Theorem \ref{main} is long and complicated. In order to
facilitate its analysis our aim in the present section is to survey
the key ideas. Theorem~\ref{main} involves a lower bound for the
number $N_X(P)$ of rational points of height at most $P$ on a quartic
hypersurface $X \subset \bfP_{\Q}^{n-1}$.  This will be achieved by
establishing an asymptotic formula for the quantity
$$
N_\omega(F;P):= \sum_{\colt{\x=(x_1,\ldots,x_n)\in \Z^n}{F(\x)=0}} 
\omega(\x/P),  
$$
as $P\rightarrow \infty$, for a suitably chosen function $\omega:\R^n
\rightarrow \R_{\geq 0}$ with compact support. 
Here  $F \in \Z[x_1,\ldots,x_n]$ denotes the quartic form
that defines $X$. 
In estimating $N_\omega(F;P)$ we will be able to recycle part 
of Birch's original argument, although there will be
a number of substantial differences. 

The starting point for the
activation of the circle method is the basic identity
\begin{equation}
  \label{eq:start}
N_\omega(F;P)=\int_0^1 S(\al) \d\al,  
\end{equation}
where $S(\al)$ is the weighted generating function
\begin{equation}
  \label{eq:S}
S(\al):=\sum_{\x\in\Z^n}\omega(\x/P) e^{2\pi i \al F(\x)},  
\end{equation}
for any $\al\in \R$. The idea is then to divide the interval $[0,1]$ into a
set of major arcs and minor arcs. Our treatment of the major arcs
follows standard procedure, and will be much in spirit with the
original argument of Birch.  It is in the treatment of the minor arcs
that our approach diverges.

Let us suppose that
\begin{equation}
  \label{eq:F}
F(x_1,\ldots,x_n)=\sum_{i,j,k,\ell=1}^n f_{ijk\ell}x_ix_jx_k x_{\ell},  
\end{equation}
for integer coefficients $f_{ijk\ell}$ that are symmetric in the indices
$i,j,k,\ell$.  Then we may define the trilinear
forms 
\begin{equation}
  \label{eq:tri}
L_i(\w;\x;\y):=4!\sum_{j,k,\ell=1}^n f_{ijk\ell}w_jx_k y_{\ell},  
\end{equation}
for $1\le i\le n$.  Using three successive applications of Weyl
differencing,  Birch ultimately relates the size of the exponential sum
$S(\al)$ to the locus of integral points on the affine variety cut out
by the system of equations $L_i(\w;\x;\y)=0$, for $1\leq i\leq n$.
This approach is quite wasteful, a fact that we are able to capitalise
on. We will use a differencing argument only once, based instead on 
the van der
Corput method, in order to relate the size of $S(\al)$ to the size of
a certain family of cubic exponential sums. 

Let $\al\in\R$, let $H\in [1,P]\cap\Z$ and write,
temporarily, 
$$
f(\x)=
\omega(\x/P)e(\alpha F(\x)).  
$$
Then the kernel of the van der Corput
method is the observation that 
$$
\#\mcal{H}S(\alpha)=\sum_{\h\in \mcal{H}}\sum_{\x\in\Z^n}f(\x+\h)=
\sum_{\x\in\Z^n}\sum_{\h\in \mcal{H}}f(\x+\h),
$$
where $\mcal{H}$ is the set of $\h\in \N^n$ such that $0<h_i\leq  H$
for $1\leq i \leq n$.
An application of Cauchy's inequality yields
\begin{align*}
H^{2n}|S(\alpha)|^2
&\ll P^n \sum_{\x\in\Z^n}\Big|\sum_{\h\in \mcal{H}}f(\x+\h)\Big|^2\\
&= P^n \sum_{\h_1\in \mcal{H}}\sum_{\h_2\in \mcal{H}}\sum_{\x\in\Z^n}
f(\x+\h_1)\overline{f(\x+\h_2)}\\
&= P^n\sum_{\h_1\in \mcal{H}}\sum_{\h_2\in \mcal{H}}\sum_{\y\in\Z^n}
f(\y+\h_1-\h_2)\overline{f(\y)}\\
&= P^n\sum_{\colt{\h\in \Z^n}{|\h|\leq H}}N(\h)\sum_{\y\in\Z^n}
f(\y+\h)\overline{f(\y)},
\end{align*}
where, in the final line, 
$$
N(\h):=\#\{\h_1, \h_2\in \mcal{H}: \h=\h_1-\h_2\} \ll H^n.
$$
We therefore conclude that
\begin{equation}\label{vdC1}
|S(\alpha)|^2\ll H^{-n}P^n\sum_{\h}|T_\h(\alpha)|
\ll \frac{P^{2n}}{H^n}+\frac{P^n}{H^{n}}\sum_{\h\neq \ma{0}}|T_\h(\alpha)|,
\end{equation}
where
$$
T_\h(\alpha):=\sum_{\x\in\Z^n}\omega_{\h}(\x/P)e\big(\al (F(\x+\h)-F(\x))\big),
$$
and 
\begin{equation}\label{wh}
\omega_\h(\x):=\omega\big(\x+P^{-1}\h\big)\omega(\x). 
\end{equation}
The reader should note that the special case $H=P$ of van~der~Corput's
method reduces to the first step in 
Birch's approach.

For each non-zero $\h \in \Z^n$ the exponential sum $T_\h(\alpha)$ is a cubic
exponential sum, involving the cubic polynomial $F_{\h}(\x)=F(\x+\h)-F(\x)$.
Note that the cubic part of $F_{\h}(\x)$ is equal to $\h.\nabla F(\x)$.
The idea will then be to estimate these exponential sums directly,
rather than using repeated applications of Weyl differencing to reduce the
degree still further. Suppose that $\al=a/q+z$ for suitable coprime integers
$a,q$.  
An application of the Poisson
summation formula will lead us to the consideration of certain
complete exponential sums modulo $q$.
These will all take the shape
$$
T_\h(a,q;\v)=\sum_{\y\bmod{q}}e_q(aF_{\h}(\y)+\v.\y),
$$
for $\v\in \Z^n$ restricted to some bounded region that expands with $P$.  
Each sum will satisfy a basic
multiplicativity property that renders it sufficient to study the sums
for prime power moduli $p^j$, for each $p^j \| q$. 

When $j=1$ or $2$ we plan to use the fact that excellent bounds exist
for $T_\h(a,p^j;\v)$ provided that the singular locus of 
the projective hypersurface $\h.\nabla
F(\x)=0$ is not too large. 
Thus we are led to make a careful study of how frequently a choice of $\h
\in \Z^n$ arises for which this singular locus has large dimension,
both as a variety over $\cQ$  and as a variety over $\overline{\F}_p$, for each $p \mid q$.
The underlying geometry of this problem will be discussed in 
\S \ref{section:geometry}.

When $j>2$ we can no longer get satisfactory individual estimates for
$T_\h(a,p^j;\v)$, even when the cubic part of $F_\h$ defines a
non-singular hypersurface over $\cQ$. 
Instead we use a more elementary argument, which bounds an average of 
sums $T_\h(a,q;\v)$, taken over a range of values for $\v$.  This 
presents considerable technical difficulties.

The necessary treatment of cubic exponential
sums will take place in \S\S \ref{s:cubic} and \ref{s:cubic_ib}. In
fact we shall treat a rather general class of cubic polynomials,
rather than restricting attention to those that arise in $T_\h(\alpha)$. This may be of independent interest. 
The remainder of the proof of Theorem \ref{main} is now relatively
straightforward. In \S \ref{s:weyl} we will review Birch's argument on
the minor arcs, in order to obtain a supplementary bound for
$S(\al)$. In \S \ref{s:vdc} we will record the bound for $S(\al)$ that
follows from the argument outlined above.  Finally in 
\S\S \ref{s:act}--\ref{s:major} we will draw to a close the proof of Theorem \ref{main}.

\section{Preliminary results}\label{section:geometry}

An important component of our work consists in viewing various
varieties that are defined over $\mathbb{Q}$ over several different
finite fields. Let us take a moment to explain precisely what we mean
by this. In general we shall be working with algebraic varieties
$W\subset \mathbb{P}_{\mathbb{Q}}^{n-1}$ that are defined by
systems of homogeneous polynomials with coefficients in
$\mathbb{Z}$. Let $\mathcal{W}=W\times_{\mathbb{Q}}\mathbb{Z}$
denote the scheme-theoretic closure of $W$ in
$\mathbb{P}_{\mathbb{Z}}^{n-1}$.  Then for any prime $p$ we can
consider the variety
$$
W_p:=\mathcal{W}\times_{\mathbb{\Z}}\mathbb{F}_p \subseteq
\mathbb{P}_{\mathbb{F}_p}^{n-1}
$$
that is obtained by reducing the coefficients of the forms defining $W$ modulo
$p$. By an abuse of notation we will normally just write $W$ for
$\mathcal{W}$ and $W_p$, it usually being clear from context which ring it
should be viewed as being defined over.
Several of the results contained in this section involve bounding the
degrees and dimensions of various varieties. 
Given an arbitrary variety $W\subset \mathbb{P}_{\mathbb{Q}}^{n-1}$,
with irreducible components $W_1,\ldots,W_D$, say, we will henceforth write
$$
\deg W := \sum_{i=1}^D \deg W_i, \quad \dim W:=\max_{1\leq i\leq
  D}\dim W_i,
$$
for the degree and dimension of $W$, respectively.
We have analogous definitions for varieties defined over finite
fields.

Let $G\in \Z[x_1,\ldots,x_n]$ be a non-zero form of degree $d\geq 2$, 
and let $p$ be a prime.  
We will find it convenient to adopt the
notation $\F_\infty:=\Q$ throughout our work.
In what follows, $v$ will always denote one of the symbols $\infty$
or $p$.  We denote by 
$$
X_G\subset\bfP_{\F_v}^{n-1}:=\proj(\F_v[x_1,\ldots,x_n]),
$$ 
the hypersurface $G=0$, viewed over $\F_v$. 
Note that when $p$ divides all of the coefficients of $G$ one will
have $X_G\cong \bfP_{\F_p}^{n-1}$ over $\F_p$. 
It will be convenient to follow the convention that the
singular locus of $X_G$ has dimension $n-1$ over
$\F_p$ when $G$ vanishes identically modulo $p$. 

Given a vector 
$$
\ma{m}=(m_1,\ldots,m_n)\in \F_v^n,
$$ 
we will write $H_{\ma{m}}\subset
\bfP_{\F_v}^{n-1}$ for the hyperplane $\ma{m}.\x=0$.  
In general, let $\sing_{\F_v}(W)$ denote the singular locus of 
a variety $W\subset \bfP_{\F_v}^{n-1}$. When $W=X_G$, we will sometimes
write $\sing_{\F_v}(G)$ for $\sing_{\F_v}(X_G)$. 
It will be convenient to set 
$$
s_v(W)=\dim \sing_{\F_v}(W),
$$
with the same convention that $s_v(G)=s_v(X_G)$.  
In particular we have $s_v(G)\in [-1,n-1]\cap \Z$, with 
$s_v(G)=n-1$ if and only if $v=p$ is such that $G$ vanishes
identically modulo $p$.
Typically we shall be interested in values of $n \geq 2$, but we shall
follow the protocol that $s_v(G)=-1$ whenever $G$ is a form in only
$1$ variable. With this in mind it is easily checked that all
of the results in this section continue to hold when $n=1$.

There is a general principle in algebraic geometry that the
basic geometric properties of $W$ over $\overline{\F}_\infty =\cQ$
should coincide with its basic properties when viewed over $\overline{\F}_p$, for $p$ sufficiently large.
In particular it is clear that $s_p(G)=s_{\infty}(G)$ for $p\gg_G 1$,
for any non-zero form $G\in \mathbb{Z}[x_1,\ldots,x_n]$, whereas one only has
$s_p(G)\geq s_{\infty}(G)$ if the prime $p$ is allowed to be arbitrary.

Given a vector $\h\in \Z^n$, the outline in  \S \ref{section:over}
gives ample evidence for the fact that we will need to control how often the form
$
\h.\nabla G(\x)
$
produces a hypersurface in $\bfP_{\F_v}^{n-1}$ with singular locus of 
prescribed dimension.
Let us write 
\begin{equation}
  \label{eq:hessian}
\ma{H}_G(\x):=\Big\{\frac{\partial^2 G}{\partial x_i
\partial x_j}\Big\}_{1\leq i,j\leq n}  
\end{equation}
for the Hessian matrix formed from the
second order partial derivatives of $G$. We observe that if $G$ is of
degree 3 then 
\begin{equation}\label{2.1a}
\ma{H}_G(\x)\y=\ma{H}_G(\y)\x. 
\end{equation}
For any $\h\in \Z^n$, we define
$$
  A_{\h}:=\{\x\in \A_{\F_v}^n: \ma{H}_G(\x)\h=\ma{0}\},
$$
where $\A_{\F_v}^{n}:=\spec(\F_v[x_1,\ldots,x_n])$ is the spectrum of 
$\F_v[x_1,\ldots,x_n]$. 
It is clear that $A_\h$ is the affine cone above 
$\sing_{\F_v}(\h.\nabla G)$. In particular $A_{\h}$ is an 
affine variety. Given any integer $s\geq 0$, we proceed to define
$$
  B_{s}:=\{\h\in \A_{\F_v}^n: \dim A_\h\geq s\}.
$$
Both of the sets $A_{\h}$ and $B_{s}$ have already been investigated
by the second author \cite[Lemma 2]{hb-india}, and the following result is a
natural generalisation of this work.

\begin{lem}\lab{lem:dimB}
Let $s$ be a non-negative integer and let $v=\infty$ or $p$, for
a prime $p\nmid d$. Then $B_s$ is an affine variety 
of degree $O_{d}(1)$, with
$$
\dim B_s \leq \min\{n, n-s+s_v(G)+1\}.
$$
\end{lem}

\begin{proof}
That $B_s$ is an affine algebraic variety of degree $O_{d}(1)$
follows immediately from the proof of \cite[Lemma
2]{hb-india}. Moreover, the upper bound $\dim B_s \leq n$ is
trivial. We may therefore proceed under the assumption that $s_v(G)<n-1$, so 
that $G$ does not vanish identically over $\F_v$. 
To obtain a second estimate for the dimension of $B_s$, we 
first show that 
\begin{equation}
  \label{eq:diag}  
\dim \{(\x,\y)\in \A_{\F_v}^{2n}: \ma{H}_G(\x)\y=\ma{0}\} \leq
n+s_v(G)+1.  
\end{equation}
An inspection of the proof of 
\cite[Lemma 2]{hb-india} reveals that this will be enough to complete
the proof of Lemma \ref{lem:dimB}. To establish \eqref{eq:diag}, we let
$S$ denote the relevant algebraic variety, and consider the diagonal 
$$
D=\{(\x,\y)\in \A_{\F_v}^{2n}: ~ \x=\y\}.
$$  
It is clear that $D$ has dimension $n$ and $D\cap S$ consists of 
all points $(\x,\x)$ for which $\nabla G(\x)=\ma{0}$, by Euler's 
identity and the fact that $p\nmid d$ if $v=p$. 
Thus it follows that $D\cap S$ has affine
dimension at most $s_v(G)+1$ in $\A_{\F_v}^{2n}$, whence 
$$
s_v(G)+1 \geq \dim D\cap S \geq \dim D+ \dim S-2n=\dim S -n,
$$
by the affine dimension theorem.  This completes the proof of \eqref{eq:diag},
and so the proof of the lemma.
\end{proof}

Recall the definition \eqref{eq:hessian} of the Hessian matrix
$\ma{H}_G$ associated to any form $G\in
\Z[x_1,\ldots,x_n]$. Then we will also need to control the size of the set 
$$
T_r:=\{\x\in \A_{\F_v}^{n}: \rank \ma{H}_{G}(\x) \leq r\},
$$
for any non-negative integer $r$.  The following result is an easy
generalisation of \cite[Lemma 2]{hb-10}, although we shall actually only employ
it in the case $v=\infty$ and $s_{\infty}(G)=-1$.

\begin{lem}\lab{lem:homrank}
Let $v=\infty$ or $p$, for
a prime $p\nmid d$, and let $r\in \Z_{\geq 0}$. 
Then $T_r$ is an affine variety of degree $O_{d}(1)$, with 
$$
\dim T_r\leq r+s_v(G)+1.
$$ 
\end{lem}

\begin{proof}
The first two claims are clear, and the third one is trivial if $G$ 
vanishes identically over $\F_v$. In order to bound the dimension of 
$T_r$ when $G$ doesn't vanish identically, we will argue by 
induction on $r\geq 0$. When $r=0$, the result is 
obvious, since then $T_0$ is contained in the set of points $\x\in
\A_{\F_v}^{n}$ for which $\nabla G(\x)=\ma{0}$, since $p\nmid d$. 
To handle the case $r\geq
1$, we will employ \eqref{eq:diag}.  To complete the proof of Lemma
\ref{lem:homrank} it will suffice to show that
\begin{equation}
  \label{eq:Urank}
  \dim U\leq r+s_v(G)+1,
\end{equation}
for any irreducible Zariski open subset $U\subseteq T_r$ on which 
$\rank \ma{H}_{G}(\x) = r$. That we may suppose the existence of such
a subset is a simple consequence of the induction hypothesis.
We now consider the incidence correspondence 
$$
I_r= \{(\x,\y)\in U\times \A_{\F_v}^{n}: \ma{H}_G(\x)\y=\ma{0}\}.
$$
The projection onto the first factor is surjective, with generic 
fibres of dimension
$n-r$. Hence it follows from \eqref{eq:diag} that
$$
\dim U=\dim I_r - (n-r)\leq  n+s_v(G)+1 -(n-r) = r+s_v(G)+1,
$$
which therefore completes the proof of \eqref{eq:Urank}.
\end{proof}

Much of our work involves
polynomials that are not necessarily homogeneous. Let $f\in 
\Z[x_1,\ldots,x_n]$ be a polynomial of degree
$d$, and write $f_0$ for the degree $d$ homogeneous part of $f$.
The following result generalises \cite[Lemma 1]{hb-10}.

\begin{lem}\lab{lem:affhess}
Let $v=\infty$ or $p$, for
a prime $p\nmid d$, and define 
$$
I:=\{(\x,\y)\in \A_{\F_v}^{2n}: \nabla^2 f(\x)\y=\ma{0}\}.
$$
Then $I$ is an affine variety of degree $O_{d}(1)$, with 
$$
\dim I\leq n+s_v(f_0)+1.
$$ 
\end{lem}

\begin{proof}
As in the previous lemma, the first two claims are obvious. To bound 
the dimension of $I$
we suppose that $f=f_0+f_1$, for a polynomial $f_1$ of degree at most $d-1$.
Then we may homogenize $f$, by writing 
\begin{equation}
  \label{eq:homogenize}
F(\x,z)=z^d f(\x/z)=f_0(\x)+z^df_1(\x/z)=
f_0(\x)+zF_1(\x,z),
\end{equation}
say, for a form $F_1(\x,z)$ of degree $d-1$.  Let
$$
J=\{(\x,\y,z)\in \bfP_{\F_v}^{2n}: \ma{H}_{f_0}(\x)\y+z
\big\{\frac{\partial^2 F_1(\x,z)}{\partial x_i\partial x_j}\big\}\y
=\ma{0}\}.
$$
Since this is merely the projective version of $I$ we see that $\dim
I=\dim J$.  On intersecting  $J$ with the hyperplane $z=0$, we deduce from
\eqref{eq:diag} that $\dim J \leq n+s_v(f_0)+1$, and Lemma 
\ref{lem:affhess} follows.
\end{proof}

In order to employ these results in our work, we will need
to be able to estimate the number of integral points on affine 
varieties that are
constrained in certain ways.  Let $W \subseteq \A_\Q^{n}$ be an 
affine variety of degree $\del$ and
dimension $\mu$. Then we have the simple upper bound
 \begin{equation}
   \label{eq:aff_triv}
\#\{\ma{t} \in W\cap\Z^{n}: |\ma{t}|\leq T\} =O_{\del}(T^\mu),
 \end{equation}
for any $T\geq 1$. This is established in \cite[Lemma 3.1]{birch}, for example.
Given a prime $p$ and an affine variety
$V\subseteq \A_{\F_p}^{n}$,  we will write $[\t]_p\in V$ to mean that the
reduction modulo $p$ of a point $\t\in \Z^n$ is contained in $V$.  
Suppose that $\Pi=\{p_1,\ldots,p_r\}$ is a finite collection of
distinct primes, with the convention that $\Pi$ is empty if $r=0$.  
Suppose further that an affine variety $V_i\subseteq
\A_{\F_{p_i}}^{n}$ is associated
to each prime $p_i\in \Pi$.   In our work we will have cause to 
estimate the quantity
$$
N(T;W,\Pi):=\#\big\{\t\in W\cap \Z^{n}: |\t|\leq T, ~\mbox{$[\t]_{p_i}\in
  V_i$ for $1\leq i\leq r$} \big\},
$$
for certain values of $T$ and certain varieties $W\subseteq
\A_{\Q}^{n}$. The outcome is the following result.
\begin{lem}\label{lem:non-trivial}
Suppose that $\dim W\leq \ell$ and $\dim V_i\leq k_i$, for $1\leq i\leq
r$, with 
$$
k_1\geq k_2\geq \cdots \geq k_r.
$$
Let $T\geq 1$ and write $D=\max\{\deg W,\deg V_i\}$. Then there exists
a constant $A(D,n)>0$ such that
$$
N(T;W,\Pi)\leq A(D,n)^{r+1}\Big(T^\ell \prod_{i=1}^r p_i^{k_i-\ell}
+ \sum_{i=1}^r T^{k_i}\prod_{j=i}^r p_j^{k_j-k_i}\Big).
$$
\end{lem}

\begin{proof} 
We prove this result by induction on the cardinality
$r=\#\Pi$. Suppose first that $r=0$. Then it follows from
\eqref{eq:aff_triv} that
$$
N(T;W,\Pi)=N(T;W,\emptyset)\ll_{D} T^{\ell},
$$ 
which is satisfactory for
the lemma.  Suppose now that $r\geq 1$.
We have two cases to consider: either $T>p_r$ or else $T\leq p_r$.

Let us deal with the case $T>p_r$ first. Our approach will be to fix a point
$\u\in V_{r}(\F_{p_r})$, and then to estimate the number of $\t$ counted by
$N(T;W,\Pi)$ for which $\t\equiv\u \bmod{p_r}$.   Now if $\t=\u+p_r
\t'$, with $|\u|\leq p_r$, then $p_r|\t'|\leq |\t-\u|\leq
T+p_r\leq 2T$. Hence $N(T;W,\Pi)$ is
\begin{align*}
&\leq \sum_{\u\in V_r(\F_{p_r})} 
\hspace{-0.1cm}
\#
\big\{\t'\in W'\cap \Z^{n}: |\t'|\leq 2T/p_r, ~\mbox{$[\t']_{p_i}\in
  V_i'$ for $1\leq i\leq r-1$} \big\}\\
&= \sum_{\u\in V_r(\F_{p_r})} N(2T/p_r;W',\Pi'),
\end{align*}
where $\Pi'=\{p_1,\ldots,p_{r-1}\}$, and the varieties $W',V_1',\ldots,
V_{r-1}'$ are obtained from $W,V_1,\ldots, V_{r-1}$ via the 
transformation $\u+p_r\t' \mapsto \t'$.
It is clear that the dimensions are 
preserved under this transformation.
Moreover, there exists a constant $A(D,n)>0$ such that
$$
\#V(\F_{p})\leq 2^{-n}A(D,n) p^{\dim V}, 
$$
for any prime $p$ and any variety $V\subseteq \bfP_{\F_p}^{n-1}$ of
degree at most $D$. Applying this with $V=V_r$ and $p=p_r$, the 
induction hypothesis therefore yields
that $N(T;W,\Pi)$ is
\begin{align*}
&\leq
\#V_r(\F_{p_r}) A(D,n)^{r}2^n 
\Big( \Big(\frac{T}{p_r}\Big)^\ell \prod_{i=1}^{r-1} p_i^{k_i-\ell}
+ \sum_{i=1}^{r-1} \Big(\frac{T}{p_r}\Big)^{k_i}\prod_{j=i}^{r-1} p_j^{k_j-k_i}
\Big)\\
&\leq A(D,n)^{r+1} \Big(T^\ell\prod_{i=1}^r p_i^{k_i-\ell} + 
\sum_{i=1}^r T^{k_i}\prod_{j=i}^r p_j^{k_j-k_i}\Big),
\end{align*}
which is satisfactory for the lemma.

Suppose now that $T\leq p_r$. We will show that 
\begin{equation}\lab{A-ind}
N(T;W,\Pi)\ll_{D}
T^{k_r}.
\end{equation}
This  will suffice to complete the proof of the lemma, since 
\begin{align*}
T^{k_r} 
&\leq  T^{k_r}+ T^\ell\prod_{i=1}^r p_i^{k_i-\ell} + \sum_{i=1}^{r-1}
T^{k_i}\prod_{j=i}^r p_j^{k_j-k_i}\\
&= T^\ell\prod_{i=1}^r p_i^{k_i-\ell} + \sum_{i=1}^{r} 
T^{k_i}\prod_{j=i}^r p_j^{k_j-k_i}.
\end{align*}
In view of the fact that $V_r$ contains at most $D$ irreducible 
components, it will
suffice to proceed under the assumption that $V_r$ is irreducible.
Our first step is to observe that 
$
N(T;W,\Pi)\leq M_{p_r,V_r}(T),
$
where $M_{p_r,V_r}(T)$ denotes the number of $\t\in \Z^{n}$ for which
$|\t|\leq T$ and $[\t]_{p_r}\in V_r$.  In order to establish 
\eqref{A-ind}, it will therefore suffice to show that 
\begin{equation}\lab{A-ind'}
M_{p,V}(T)\ll_{D}T^{\dim V},
\end{equation}
for any prime $p$, and any irreducible affine variety
$V\subseteq \A_{\F_p}^{n}$ of degree $D$. We will establish
\eqref{A-ind'} by induction on the dimension $\mu$ of $V$.  
Since an irreducible variety of dimension zero contains just one 
point, the estimate is trivial when $\mu=0$.  
Assume now that $\mu\geq 1$.  Since $V$ is assumed to be irreducible
we may find an index $1 \leq a\leq n$ such 
that $V$ intersects the hyperplane $t_a=\al$ properly, for every 
$\al$. Let $H_\al$ denote this
hyperplane.  In this way we obtain the upper bound
$$
M_{p,V}(T) \leq \sum_{|\al| \leq T} M_{p, V\cap H_\al}(T).
$$
Since $V\cap H_\al$ has dimension at most $\mu-1$ for every $\al$, 
and decomposes into at most $D$ irreducible components, 
an application of the  induction hypothesis implies that  
$M_{p, V\cap H_\al}(T)=O_{D}(T^{\mu-1})$.  This suffices to complete 
the proof of
\eqref{A-ind'}, and so completes the proof of the lemma.
\end{proof}

Taking $W=\A_{\Q}^n$ and $T=p_1\cdots p_r$ in the statement of Lemma
\ref{lem:non-trivial}, it is now a trivial matter to deduce that there
is a constant $A(D,n)>0$ such that
\begin{equation}
  \label{eq:n-tcor}
\#\big\{\t \bmod{p_1\cdots p_r}: \mbox{$[\t]_{p_i}\in V_i$ 
for $1\leq i\leq r$} \big\}
\leq A(D,n)^r\prod_{i=1}^r p_i^{\dim V_i},
\end{equation}
where $D=\max \{\deg V_i\}$. 

Several of the arguments that we will encounter involve inducting on
the dimension of the singular locus of certain varieties. 
The following result will prove extremely useful in this context.

\begin{lem}\lab{lem:hyp}
Let $\Pi$ be a collection of primes, with $\#\Pi=r\geq 0$, 
and write $\Pi_a:=\{p\in\Pi: p> a\}$ for each $a\in\N$.
Then there exists a primitive vector
$\ma{m}\in\Z^n$ and a constant $c=c(d,n)>0$ such
that the following hold:
\begin{enumerate}
\item $\ma{m}\ll_{d} r$, 
\item 
for each $v\in
  \{\infty\}\cup \Pi_{cr}$, we have
  \begin{equation}
    \label{eq:hyp'}
s_{v}(X_G\cap H_{\ma{m}})=
\left\{
\begin{array}{ll}
-1, & \mbox{if $s_v(X_G)=-1$,}\\ 
s_{v}(X_G)-1,
& \mbox{if $s_v(X_G)\geq 0$,} 
\end{array}
\right.
  \end{equation}
\item
for any $\ma{e}\in \Z^n\setminus \{\ma{0}\}$ such that 
$\ma{m}.\ma{e}=0$, we have $|\ma{e}|\gg |\ma{m}|^{1/(n-1)}$.
\end{enumerate}
\end{lem}

\begin{proof}
When $r=0$ this is just \cite[Lemma 4]{hb-india}. We proceed
under the assumption that $r\geq 1$, so that $\Pi$ is non-empty. 
Let us begin by showing that 
\begin{equation}\lab{lower-s}
s_{v}(X_G\cap H)\geq \max\{s_{v}(X_G)-1,-1\},
\end{equation}
for any $v\in \{\infty\}\cup \Pi$ and any hyperplane $H$.  
Now it  is easy to check that
$\sing_{\F_v} (X_G\cap H)$ consists of all points $[\x]\in X_G\cap
H$ for which $\nabla G(\x)$ is proportional to the
coefficient vector defining $H$.  In particular we have 
$$
\sing_{\F_v} (X_G) \cap H \subseteq \sing_{\F_v} (X_G\cap H).
$$
Thus it follows that 
$$
s_{v}(X_G\cap H)\geq  \dim \big(\sing_{\F_v}(X_G)\cap H\big) \geq  
\max\{s_{v}(X_G)-1-1\}. 
$$
Here the lower bound $s_{v}(X_G)-1$ follows from 
the projective dimension theorem and the lower bound $-1$ is trivial.  
This therefore establishes \eqref{lower-s}.

In what follows, let us write $Z^n$ for the set of primitive
vectors in $\Z^n$. Let $\ma{m}\in Z^{n}$ and let $H$ be the
hyperplane $\ma{m}.\x=0$.  Let $v\in \{\infty\}\cup\Pi_d$. 
We will derive the existence of a proper subvariety 
$W\subset  \bfP_{\F_v}^{n-1}$ of degree $O_{d}(1)$, such that $[\ma{m}]_v\in
W(\F_v)$ whenever there is strict inequality in \eqref{lower-s}.
It is clear from the above that this can only happen if 
\begin{equation}
  \label{eq:eq}
\dim \big(\sing_{\F_v}(X_G)\cap H\big) = s_{v}(X_G)\geq 0,  
\end{equation}
or if 
\begin{equation}
  \label{eq:eq'}
s_{v}(X_G\cap H)>\dim \big(\sing_{\F_v}(X_G)\cap H\big).
\end{equation}
Let $\widehat{X_G}\subset \bfP_{\F_v}^{n-1}$ denote
the dual variety, which will be a proper subvariety of degree 
$O_{d}(1)$, satisfying the property that 
$[\y]_v\in \widehat{X_G}$ if there exists $[\ma{z}]_v\in X_G$ 
such that $\y=\nabla G(\ma{z})$. 
However, for \eqref{eq:eq'} to occur  we would need 
there to exist a point $[\ma{z}]_v\in \sing_{\F_v}(X_G\cap H)$,
which is not contained in the singular locus of $X_G$. But this
means that $\nabla G(\ma{z})$ must be a non-zero multiple of $\ma{m}$, whence 
$[\ma{m}]_v\in \widehat{X_G}(\F_v)$.

Suppose now that \eqref{eq:eq} holds. Then it 
follows that $H$ contains an irreducible
component of $\sing_{\F_v}(X_G)$ of maximal dimension.  
Suppose that this singular locus has $D$ such components. 
Then $D=O_{d}(1)$, since the degree of $\sing_{\F_v}(X_G)$ can be
bounded in terms of $d$ and $n$.
On picking points 
$$
[\y_1]_v,\ldots,[\y_D]_v\in \sing_{\F_v}(X_G),
$$ 
one from each
component, we see that $[\ma{m}]_v$ will be contained in the union of
hyperplanes $\ma{y}_i.\x=0$, for $1\leq i\leq D$.

In the case that $s_v(X_G)\geq 0$, we will take $W\subset
\bfP_{\F_v}^{n-1}$ to be the union of
$\widehat{X_G}$ and all these hyperplanes. When  
$s_v(X_G)=-1$, we will take $W=\widehat{X_G}$. 
In conclusion, we have constructed a proper subvariety
$W\subset \bfP_{\F_v}^{n-1}$ of degree $O_{d}(1)$, such that $[\ma{m}]_v\in
W(\F_v)$ whenever 
\eqref{eq:hyp'} is false.
Our argument produces such a variety $W=W_v$, say, for each $v\in
\{\infty\}\cup\Pi_d.$

Thus far our work will ultimately allow us to choose a vector $\ma{m}\in
Z^n$, with low height, such that the second condition is
satisfied in the statement of the lemma. In order to handle the final
condition, we introduce the set
$$
\{\ma{m}\in Z^n: |\ma{m}|\leq M, ~\mbox{$\exists ~\ma{e}\in Z^n$
  such that $\ma{m}.\ma{e}=0$ and $|\ma{e}|\leq A |\ma{m}|^{1/(n-1)}$}\},
$$
for given $A>0$ and $M\geq 1$. Let us denote this set by $S_A(M)$.
We claim that there exists a constant $C_1>0$ depending only on $n$, such that 
\begin{equation}
  \label{eq:S_upper}
  \#S_A(M)\leq A^{n-1}C_1 M^n.
\end{equation}
This is trivial for $A\ge 1$, so we shall assume that $A<1$.
To prove \eqref{eq:S_upper} for $A<1$ we can apply the 
geometry of numbers, and the facts that we will
need may all be read off from \cite[Lemma 1]{annal}.  Breaking the
possible ranges for $\ma{e}$ into dyadic intervals, we obtain
$$
\#S_A(M) \leq \sum_{E\leq AM^{1/(n-1)}} \sum_{\colt{\ma{e}\in
    Z^n}{E<|\ma{e}|\leq 2E}} \#\{\ma{m}\in Z^n: |\ma{m}|\leq M, 
~\ma{m}.\ma{e}=0\}, 
$$
where $E$ runs over powers of $2$. 
The vectors $\ma{m}$ appearing in the summand are restricted to an
integer lattice of rank $n-1$ and determinant $|\ma{e}|>E$. 
It follows that 
$$
\#S_A(M) \ll \sum_{E\leq AM^{1/(n-1)}} 
\hspace{-0.2cm}
\sum_{\colt{\ma{e}\in
    Z^n}{E<|\ma{e}|\leq 2E}}  
\hspace{-0.2cm}
\frac{M^{n-1}}{E}
\ll
\hspace{-0.1cm}
\sum_{E\leq AM^{1/(n-1)}} 
\hspace{-0.2cm}
(EM)^{n-1} \ll A^{n-1}M^n,
$$
which thereby completes the proof of \eqref{eq:S_upper}.

Now the number of  vectors $\ma{m}\in
Z^{n}$ such that  $|\ma{m}|\leq M$ is at least $C_2 M^{n}$ for some
constant $C_2>0$, depending only on $n$.  Let $p\in \Pi$ be any prime.
Lemma \ref{lem:non-trivial} reveals that the number of vectors $\ma{m}\in
Z^{n}$ such that  $|\ma{m}|\leq M$ and either $\ma{m}\in W_\infty$ or
$[\ma{m}]_{p}\in W_{p}(\F_{p})$, is 
\begin{equation}
  \label{eq:suit}
\leq 
C_3
\big(\frac{M^{n}}{p}+M^{n-1}\big)
\end{equation}
for some constant $C_3>0$ depending only on $d$ and $n$. 
We now focus our attention on the primes that are contained in the set
$\Pi_{3rC_3/C_2}=\Pi'$, say. Recall that these are the primes in $\Pi$ that
exceed $3rC_3/C_2$ in size, so that in particular,
$$
\sum_{p\in \Pi'}\frac{1}{p}< \frac{rC_2}{3rC_3}=\frac{C_2}{3C_3}.
$$
Set $A:=(C_2/3C_1)^{1/(n-1)}$.  Then it follows from combining 
\eqref{eq:S_upper} with \eqref{eq:suit} that
the number of vectors $\ma{m}\in Z^{n}$ such that  $|\ma{m}|\leq M$
and either $\ma{m}\in W_\infty\cup S_{A}(M)$, or else
$[\ma{m}]_{p}\in W_{p}(\F_{p})$ for some $p\in\Pi'$, is 
\begin{align*}
\leq 
A^{n-1}C_1M^n+
C_3\sum_{p\in \Pi'}\big(\frac{M^{n}}{p}+M^{n-1}\big)
&< C_3rM^{n-1}+\frac{2C_2M^n}{3}.
\end{align*}
It is now clear that by taking $M$ to be any integer $\geq 3C_3r/C_2$,
we can ensure that this upper bound is strictly less than $C_2M^n$. But
this means that there exists a primitive vector $\ma{m}\in\Z^{n}$, with
$\ma{m}\ll_{d} r$, such that  
\eqref{eq:hyp'} holds for each $v\in\{\infty\}\cup\Pi'$, and also 
$|\ma{e}|\gg |\ma{m}|^{1/(n-1)}$ for any $\ma{e}\in Z^n$ such that
$\ma{m}.\ma{e}=0$.  This completes the proof of Lemma \ref{lem:hyp}.
\end{proof}

The essential content of Lemma \ref{lem:hyp} is that we can always
find a primitive vector $\ma{m}\in \Z^n$, of  low height,
such that the dimension of the singular locus of $X_G\cap
H_{\ma{m}}$ is simultaneously reduced by $1$ 
over many ground fields, at least if $X_G$ is singular.
The final condition on $\ma{m}$ appearing in
the lemma might seem out of place, but its presence affords us better
control over the effect that taking hyperplane sections has on 
certain weight functions used in our work.

We end this section with some basic estimates for exponential sums
and integrals over arbitrary polynomials.  Let $f\in \Z[x_1,\ldots,x_n]$ be a
polynomial of degree $d$, and let $p$ be a prime. 
As previously we write $f_0$ for the degree $d$  
homogeneous part of $f$,
and $s_p(f_0)$ for the dimension of the
singular locus of the hypersurface $f_0=0$ in $\bfP_{\F_p}^{n-1}$.
Let $\|f\|$ denote the maximum modulus of the coefficients of $f$.

We begin by looking at a certain family of weighted exponential integrals. 
It will be convenient to work with infinitely differentiable weight functions
$w: \R^n \rightarrow \R_{\geq 0}$, which have compact support and take
non-negative real values. Given such a function $w$, 
we set $S(w)$ to be the smallest $S$
such that $w$ is supported in the hypercube $[-S,S]^n$, and we let
$$
S_j(w):=\max\Big\{ \Big|
\frac{\partial^{j_1+\cdots+j_n}w(\x)}{\partial^{j_1}x_1\cdots
\partial^{j_n}x_n}\Big|: ~\x \in \R^n, ~j_1+\cdots+j_n=j\Big\},
$$
for each integer $j \geq 0$.  Let constants $c_n$ and $c_{n,j}$ be
given, and define $\WW_n$ to be 
the set of infinitely differentiable functions $w: \R^n \rightarrow
\R_{\geq 0}$ of compact support, such that $S(w)\le c_n$ and
$S_j(w)\le c_{n,j}$ 
for all $j\ge 0$.  In future all our order
constants will be allowed to depend on $c_n$ and the $c_{n,j}$,
without further comment.
Our work will lead us to consider integrals of the shape 
\begin{equation}
  \label{eq:Iw} 
I=I_w(f,t,\ma{u}):=\int_{\R^n}w(\x)e(t f(\x)-\ma{u}.\x)\d\x, 
\end{equation}
for given $w\in\WW_n$, $t\in\R$ and $\u\in \R^n$.
The following result is a straightforward consequence of the second author's
work \cite[Lemma 10]{HB-crelle}.

\begin{lem}\lab{lem:Iq}
Let $f\in \Z[\x]$ be a polynomial of degree $d$, and let $w\in \WW_n$
with $S(w)\le S$.
Let $H\geq \|f\|$ and $N\geq 1$.  Then there exists a constant $c(d,n)>0$
such that 
\begin{align*}
I\ll_{d,N} 
\frac{1}{R^N}+
\meas\big\{\x:~|\x|\leq 1+S, ~|t \nabla
f(\x) - \ma{u}|\leq R \max\{1,\sqrt{|t|H}\}\big\},
\end{align*}
for any $R\geq c(d,n)$.
\end{lem}

\begin{proof}
Let $\del\in \R$ be a parameter in the range $0<\del\leq 1$, to be
chosen in due course. 
To prove Lemma \ref{lem:Iq} we will decompose the integral
$I=I_{w}(f,t,\ma{u})$ into a number of integrals over subregions
of size at most $\delta$. It follows from \cite[Lemma 2]{HB-crelle}
that there exists a weight function
$w_\del(\x,\y)\in \WW_{2n}$, for suitable constants $c_{2n},
c_{2n,j}$, 
such that
$$
w(\x)=\del^{-n}\int w_{\del}\Big(\frac{\x-\y}{\del},\y\Big)\d\y.
$$
Moreover, ${\rm supp}(w_\del)\subseteq [-1,1]^n\times [-1-S,1+S]^n$.
Write $\mcal{B}=[-1-S,1+S]^n$, for ease of notation. 
Then on making this substitution into $I$, and writing $\x=\y+\del \z$, 
we obtain
\begin{align}
\label{eq:lemIq:1}
|I|
&=\del^{-n}\Big|\int\int
w_\del\big(\del^{-1}(\x-\y),\y\big)e(tf(\x)-\ma{u}.\x)\d\x\d\y\Big|\nonumber\\
&\leq\int_{\mcal{B}}\Big|\int w_\del(\z,\y)e(tf(\y+\del \z)-\del 
\ma{u}.\z)\d\z\Big|\d\y\nonumber\\
&=\int_{\mcal{B}} |J(\y)|\d\y,
\end{align}
say. 
Write $F(\z)=tf(\y+\del \z)-\del \ma{u}.\z$, for fixed $\y$, and
recall that $S=O(1)$. Then it is easy to see that the $k$th power
derivatives of $F(\z)$ are all $O_k(\del^k |t|H)$, for $k\geq 2$, when
$(\y,\z)\in \mcal{B}\times [-1,1]^n$.  
Similarly, one finds that
$$
\nabla F(\z) = \del t\nabla f(\y) - \del \ma{u} +O(\del^2 |t|H).
$$
When $|t\nabla f(\y) - \ma{u}|$ is large we will be able to show that
$J(\y)$ is small in  \eqref{eq:lemIq:1}. Alternatively, we will
simply use the trivial bound $J(\y)\ll 1$.

Let $R\geq 1$ and suppose that $\y$ is such that $|t\nabla f(\y) -
\ma{u}|\geq \del^{-1}R$.  Then it follows from our work above that 
there exists a constant $c(d,n)>0$
such that $|\nabla F(\z)|\gg R$, provided that 
$$
R\geq c(d,n) \del^2 |t| H.
$$
We will take $\del=\min\{1,(|t|H)^{-1/2}\}$, so that $0<\del\leq
1$. An application of \cite[Lemma 10]{HB-crelle} now reveals that
$J(\y)\ll_{d,N} R^{-N}$ for any $N\geq 1$, when $R\geq c(d,n)$.
We may now insert this into \eqref{eq:lemIq:1} to deduce that
$$
I \ll_{d,N} R^{-N}+ 
\meas\big\{\y\in\mcal{B}:~ |t \nabla 
f(\y) - \ma{u}|\leq R \max\{1,\sqrt{|t|H}\}\big\},
$$
for any $N\geq 1$ and $R\geq c(d,n)$. 
This completes the proof of the lemma.
\end{proof}

We will also need good upper bounds for complete exponential
sums modulo $p$ or $p^2$.

\begin{lem}\label{deligne}
Let $f\in \Z[\x]$ be a
polynomial of degree $d\geq 2$, and let $p$ be a prime. 
Then we have
$$
\sum_{\x\bmod{p^j}} e_{p^j}(f(\x)) \ll_d p^{j(n+1+s_p(f_0))/2},
$$
for $j=1,2$.
\end{lem}
\begin{proof} 
When $j=1$ this can be extracted from the work of Hooley
\cite{hooley-38}. It is established by induction on $s_p(f_0)$, the
inductive base $s_p(f_0)=-1$ being taken care of by Deligne's 
estimate. The general case
is reduced to this situation by appropriate hyperplane sections. 
The result is trivial when $p\mid d$ or if $p$ divides all of the
coefficients of $f_0$, and so we proceed under the 
assumption  that $p\nmid df_0$. 

When $j=2$ we may write $\x=\y+p \z$, giving
\begin{align*}
\sum_{\x\bmod{p^2}} e_{p^2}(f(\x)) 
&=
\sum_{\y,\z\bmod{p}} e_{p^2}(f(\y))e_p(\z.\nabla f(\y))\\
&\ll
p^n \#\{\x \bmod{p}: p\mid \nabla f(\x)\}.
\end{align*}
Suppose that $f=f_0+f_1$, for a polynomial $f_1$ of degree at most $d-1$.
Arguing as in the proof of Lemma \ref{lem:affhess}, we 
homogenize $f$, giving \eqref{eq:homogenize},
with a form $F_1(\x,z)$ of degree $d-1$.   It now follows that
\begin{align*}
\sum_{\x\bmod{p^2}} e_{p^2}(f(\x)) 
&\ll
p^n \#\Big\{(\x,1) \bmod{p}: p\mid \frac{\partial F}{\partial x_i}(\x,1),
~(1\leq i\leq n)\Big\}\\
&\ll_{d,n} p^{n+\dim V},
\end{align*}
where $V$ denotes the projective variety defined by 
$\partial F/\partial x_i(\x,z)=0$, for $1\leq i\leq n$.
To complete the proof of Lemma
\ref{deligne}, it therefore suffices to show that
$
\dim V\leq 1+s_{p}(f_0).
$  
But this follows immediately on noting that the intersection of 
$V$ with the hyperplane
$z=0$ is just $\sing_{\F_p}(f_0)$.
\end{proof}

\section{Cubic exponential sums: the main estimate}\label{s:cubic}

The focus of the paper now shifts towards estimating 
a rather general family of cubic exponential sums. 
Let $g\in \Z[x_1,\ldots,x_n]$ be an arbitrary cubic polynomial, which
is not necessarily homogeneous. The central object of study is the 
exponential sum
\begin{equation}\label{eq:Tcubic}
\TT(\al)=\TT_n(\al;g,w,P):=\sum_{\x\in\Z^n}w(\x/P)e(\al g(\x)),
\end{equation}
for a suitable family of weights $w$ on $\R^n$. Recall the definition
of the set $\WW_n$ of infinitely differentiable weight functions 
$w: \R^n \rightarrow
\R_{\geq 0}$, that were introduced in the preceding section.  
Given $P \geq 1$, we will need to
work with the function
$$
\|g\|_P:= \|P^{-3}g(Px_1, \ldots, Px_n)\|,  
$$
where $\|f\|$ denotes the usual height of a polynomial $f$. In
particular, it is clear that 
\begin{equation}
  \label{eq:cots}
\|g_0\|=\|g_0\|_P\leq \|g\|_P\leq \|g\|,
\end{equation}
for any $P\geq 1$, 
where $g_0$ denotes the cubic homogeneous part of $g$.

We are almost ready to reveal our first bound for \eqref{eq:Tcubic}.
We will assume throughout this section and the next that
$\al=a/q+z$, with $a,q\in \Z$  such that 
\begin{equation}
  \label{eq:aqP}
1\leq a\leq q\leq P^2, \quad \hcf(a,q)=1,
\end{equation}
and $z\in\R$ such that 
\begin{equation}\label{eq:z}
  |z|\leq q^{-1}P^{-1}.
\end{equation}
It will be convenient to set
$$
s_\infty:=s_\infty(g_0), \quad s_p:=s_p(g_0),
$$
for each $p\mid q$. 
We will write $q=bc^2d$, where
\begin{equation}
  \label{eq:bcd}
b:=\prod_{\colt{p^e\| q}{e\leq 2}}p^e, 
\quad
d:=\prod_{\colt{p^e\| q}{e\geq 3, ~2\nmid e}}p.
\end{equation}
It is not hard to see that $d$ divides $c$, and that 
there exist a divisor $d_0$ of $d$ such that
$d_0^{-1}d^{-1}c$ is a square-full integer. Moreover, 
$\hcf(b,c^2d)=1$. Finally, we define
\begin{equation}
  \label{eq:qi}
  r_i:=\prod_{\colt{p^e\| bd}{s_p=i-1}}p^{e},
\end{equation}
for $0\leq i\leq n$. We have the following result.

\begin{pro}\label{lem:T1}
Let $A,\ve>0$. Let $w \in \WW_n$ 
and let $g\in \Z[x_1,\ldots,x_n]$ be a cubic polynomial with 
$\|g\|_P\leq H$, for some
$H$ in the range $1\leq H\leq P^A.$ Assume that $s_\infty=-1$. 
Let $a,q$ be such that \eqref{eq:aqP} holds and 
$q=bc^2d$, in the notation of \eqref{eq:bcd}. Define
\begin{equation}
  \label{eq:V}
V:=qP^{-1} \max\{1,\sqrt{|z|HP^3}\},
\end{equation}
and
\begin{equation}
  \label{eq:W}
  W:=V+ \min\big\{(c^2d H)^{1/3}, c^{1/2}V^{1/2}+ c^{5/6}H^{1/6}\big\}.
\end{equation}
Then we have
$$
\TT(a/q+z)\ll_{A}q^{-n/2}
\Big(\prod_{i=0}^n r_{i}^{i/2}\Big)
P^{n+\ve}W^{n}.
$$
\end{pro}

The proof of \prop \ref{lem:T1} will be carried out in \S 
\ref{s:cubic_ib}. Using an argument based on induction we are now  
in a position to build on this result, in order to establish the 
following generalisation.

\begin{pro}\label{lem:T2}
Let $A,\ve>0$. Let $w \in \WW_n$  
and let $g\in \Z[x_1,\ldots,x_n]$ be a cubic polynomial with 
$\|g\|_P\leq H$, for some
$H$ in the range $1\leq H\leq P^A.$
Let $a,q$ be such that \eqref{eq:aqP} holds and 
$q=bc^2d$, in the notation of \eqref{eq:bcd}. Then we have 
$$
\TT(a/q+z)\ll_{A} \min_{1+s_\infty\leq \eta \leq n}
q^{-(n-\eta)/2}\Big(\prod_{i=\eta}^n r_{i}^{(i-\eta)/2}\Big)
P^{n+\ve}W^{n-\eta},
$$
where $W$ is given by \eqref{eq:W}.
\end{pro}

Note that Proposition \ref{lem:T2} applies to cubic polynomials with
arbitrary singular locus, whereas
Proposition \ref{lem:T1} is only valid when $g_0$ is a non-singular
cubic form. In fact the statement of Proposition \ref{lem:T1} is
retrieved by taking 
$s_{\infty}=-1$ and $\eta=0$ in Proposition \ref{lem:T2}.

In order to establish \prop \ref{lem:T2}, it will suffice to show that
$$
\TT(a/q+z)\ll_{A} 
q^{-(n-\eta)/2}
\Big(\prod_{i=\eta}^n r_{i}^{(i-\eta)/2}\Big)
P^{n+\ve}W^{n-\eta},
$$
for any integer $\eta$ in the interval $[1+s_\infty,n]$.  
The proof of this estimate will be by induction on
$n$. The case $n=1$, for which $s_\infty=-1$, is handled
by \prop \ref{lem:T1} when $\eta=0$. When $n=\eta=1$ the bound is
trivial since we always have $\TT(a/q+z)\ll P$.
We proceed under the assumption that $n \geq 2$, and our induction
hypothesis is that
\begin{equation}
  \label{eq:suffice'}
\TT_{n-1}(\al;h,w_0,P)
\ll_{A} 
q^{-(n-1-\eta)/2}
\Big(\prod_{i=\eta}^{n-1} r_{i}^{(i-\eta)/2}\Big)
P^{n-1+\ve}W^{n-1-\eta},
\end{equation}
for any integer $\eta$ in the interval $[1+s_\infty(h),n-1]$, 
and any suitable  $w_0\in
\WW_{n-1}$ and cubic polynomial $h\in \Z[x_1,\ldots,x_{n-1}]$.

Let $w\in \WW_n$ and let $g \in \Z[x_1,\ldots,x_n]$ be a cubic
polynomial. In particular, we have $s_p\geq s_{\infty}\geq -1,$
for any prime $p$.
Let $P\geq 1$ and let $H$ be such that $\|g\|_P\leq H\leq P^A$, for
some $A>0$.    Our plan will be to use hyperplane sections,
in order to reduce the problem to a
consideration of $(n-1)$-dimensional exponential sums involving 
cubic polynomials
whose cubic part defines a hypersurface with 
singular locus of dimension $\max\{-1,s_{\infty}-1\}$.
We take $\Pi$ to be the set of primes $p\mid q$, and set
$$
r:=\#\Pi=\omega(q).
$$
In particular $r\ll P^\ve$, by \eqref{eq:aqP}. According to Lemma 
\ref{lem:hyp}, there
exists a constant $c=c(n)>0$  and a primitive vector 
$\ma{m}\in\Z^n$, with $\ma{m}\ll r$,  
such that 
$$
s_{v}(X_{g_0}\cap H_{\ma{m}})= \max\{-1,s_v(g_0)-1\},
$$
for each $v \in  \{\infty\}\cup \Pi_{cr}$, and with
$|\ma{e}|\gg |\ma{m}|^{1/(n-1)}$ for any 
$\ma{e}\in \Z^n\setminus\{\ma{0}\}$ such that $\ma{m}.\ma{e}=0$. 
In order to apply the induction hypothesis we will sum over affine
hyperplane sections $\ma{m}.\x=k$, for integers $k\ll Pr\ll P^{1+\ve}$.
This gives 
\begin{align}\lab{eq:ind}
\TT_n(a/q+z;g,w,P)
= &\sum_{k\ll P^{1+\ve}} 
\sum_{\colt{\x\in\Z^n}{\ma{m}.\x=k}} w(\x/P) e(\al g(\x))\nonumber\\
=&\sum_{k\ll P^{1+\ve}} \mcal{S}_k,
\end{align}
say. Now $\mcal{S}_k$ is zero unless there exists a vector $\ma{t}\in \Z^n$
such that $\ma{m}.\ma{t}=k$ and $\ma{t}\ll P$. Let us fix such a
choice of vector, and write $\x=\ma{t}+\y$ in $\mcal{S}_k$.
Then clearly $\ma{m}.\x=k$ if and only if $\ma{m}.\y=0$. This
condition defines a lattice $\sfl\subseteq \Z^n$ of rank $n-1$ and
determinant $|\ma{m}|$, by part (i) of \cite[Lemma~1]{annal}. We now choose a
basis $\ma{e}_1,\ldots, \ma{e}_{n-1}$ for $\sfl$, 
as in part (iii) of \cite[Lemma 1]{annal}. 
Then our choice of $\ma{m}$ ensures that 
for each $1\leq i\leq n-1$, we have
$$
L\ll |\ma{e}_i| \ll L,
$$
where $L:=|\ma{m}|^{1/(n-1)}$.
Moreover, any of the vectors $\y$ we are interested in can be written as
$\y=\sum_{i=1}^{n-1}\la_i \ma{e}_i$ for 
$\mla=(\la_1,\ldots,\la_{n-1})\in\Z^{n-1}$ such that
$$
\la_i \ll  \frac{P}{|\ma{e}_i|} \ll \frac{P}{L}, \quad(1\leq i\leq n-1).
$$
It will suffice to assume that $L\ll P$ in  what follows. Indeed, the
alternative hypothesis implies that $P\ll L \ll r^{1/(n-1)}\ll P^\ve$,
which is a contradiction for large enough $P$.
Putting all of this together, we conclude that
\begin{align}
  \label{eq:ind'}
  \mcal{S}_k&= \sum_{\mla\ll P/L} 
w\Big(\big(\t+\sum_{i=1}^{n-1}\la_i \ma{e}_i\big)/P\Big)e\Big(\al
g\big(\t+\sum_{i=1}^{n-1}\la_i \ma{e}_i\big) \Big)\nonumber\\
&=  \sum_{\mla\in\Z^{n-1}} 
w_0\big(L\mla /P\big) e\big(\al h(\mla) \big)\nonumber\\
&=
\TT_{n-1}(a/q+z;h,w_0,P/L),
\end{align}
where 
\begin{equation}
  \label{eq:hw}
h(\u):=g\Big(\t+\sum_{i=1}^{n-1}u_i \ma{e}_i\Big), \quad
w_0(\ma{u}):=
w\Big(P^{-1}\t+L^{-1}\sum_{i=1}^{n-1}u_i \ma{e}_i\Big),
\end{equation}
and $\ma{u}=(u_1,\ldots,u_{n-1})$.
Our task is to show that we can apply the
induction hypothesis to estimate $\TT_{n-1}(a/q+z;h,w_0,P/L)$.

We claim that $w_0\in\WW_{n-1}$, for our choice of $\ma{m}$ and
$\ma{t}$.  Now it is clear that $w_0:\R^{n-1}\rightarrow
\R_{\geq 0}$ is an infinitely
differentiable function, such that
$S_j(w_0)\ll_j 1$ for each $j\geq 0$. Thus it remains to show that
$w_0$ has compact support, with 
$S(w_0)\ll 1$. Let $\u\in \R^{n-1}$ be a non-zero vector such that 
$w_0(\ma{u})\neq 0$. We
wish to show that $\u\ll 1$. Now it is clear that $w_0(\ma{u})=0$ unless
$$
P^{-1}\t+L^{-1}\sum_{i=1}^{n-1}u_i \ma{e}_i\ll 1,
$$
since $S(w)\ll 1$.  The bound $\t\ll P$ then implies that
$\sum_{i=1}^{n-1}u_i \ma{e}_i \ll L$.  If $v_i$ is the integer part
of $u_i$ we now have $\sum_{i=1}^{n-1}v_i \ma{e}_i\ll L$, since
$\ma{e}_i\ll L$.  It then follows from part(iii) of \cite[Lemma
8]{annal} that $v_i\ll L/|\ma{e}_i|\ll 1$, whence $u_i\ll 1$, as
required.

We now turn to the cubic polynomial $h\in\Z[u_1,\ldots,u_{n-1}]$
defined in \eqref{eq:hw}.
Note first that  
\begin{align*}
\|h\|_{P/L} &=
\|L^3P^{-3}g(\t+Pu_1\ma{e}_1/L+\cdots +Pu_{n-1}\ma{e}_{n-1}/L)\|\\ 
&\ll L^3 \|g\|_P\\
&\leq H L^3,
\end{align*}
since $\t\ll P$ and $\ma{e}_i\ll L$.
Our final task is to show that 
$$
s_{v}(h_0)=\max\{-1,s_{v}(g_0)-1\},
$$
for each $v \in  \{\infty\}\cup \Pi_{cr}$.
Define the $n\times (n-1)$ matrix $\ma{E}$ to have column vectors
$\ma{e}_1,\ldots \ma{e}_{n-1}$.  Then 
$h(\u)=g(\ma{E}\u+\t)$. It is not hard to see  that 
the homogeneous cubic part of $h(\u)$ is just
$h_0(\u)=g_0(\ma{E}\u)$.
Viewed over $\F_v$, for each $v \in  \{\infty\}\cup \Pi_{cr}$, we see
that the locus of $[\u]\in \bfP_{\F_v}^{n-2}$ such that  
$g_0(\ma{E}\u)=0$ is isomorphic to the locus of  
$[\y]\in \bfP_{\F_v}^{n-1}$ 
such that $g_0(\y)=\ma{m}.\y=0$. This
therefore establishes the claim, since we have already seen that
$s_{v}(X_{g_0}\cap H_{\ma{m}})= \max\{-1,s_v(g_0)-1\}$.

We are now in a position to apply the induction hypothesis
\eqref{eq:suffice'} to
estimate the quantity $\TT_{n-1}(a/q+z;h,w_0,P/L)$ in 
\eqref{eq:ind'}, with $H$ replaced by 
$cHL^3$ for a suitable constant
$c\ll 1$. In particular we have $cHL^3\leq P^{A+4}$, provided that $P$
is sufficiently large.  On recalling that 
$$
1\leq L=|\ma{m}|^{1/(n-1)}\ll r^{1/(n-1)}\ll P^{\ve},
$$ 
we therefore deduce from \eqref{eq:ind} and the induction hypothesis
with $\eta$ replaced by $\eta-1$ that 
\begin{equation}
  \label{eq:induce}
\TT(a/q+z)\ll_{A} 
q^{-(n-\eta)/2}P^{n+\ve}W^{n-\eta}
\prod_{i=\eta-1}^{n-1} q_{i}^{(i-\eta-1)/2},
\end{equation}
for any integer $\eta$ in the interval $[2+\max\{-1,s_\infty-1\},n]$, 
where $s_\infty=s_\infty(g_0)$, and 
$$
q_i:= \prod_{\colt{p^e\| bd}{s_p(h_0)=i-1}}p^{e}.
$$
Suppose first that $s_{\infty}\geq 0$. Then the above estimate holds
for integers $\eta \in [1+s_\infty,n]$, as required. 
Moreover we also have $s_p\geq 0$. Thus it
follows that 
\begin{align*}
q_i= 
\prod_{\tstack{p^e\| bd}{s_p(h_0)=i-1}{p\leq c r}}p^{e}
\prod_{\tstack{p^e\| bd}{s_p(h_0)=i-1}{p>cr}}p^{e}
&=
\prod_{\tstack{p^e\| bd}{s_p(h_0)=i-1}{p\leq c r}}p^{e}
\prod_{\tstack{p^e\| bd}{s_p(g_0)=i}{p>cr}}p^{e}\\
&\leq 
r_{i+1}
\prod_{p\leq c r}p^{2}\\
&\leq (cr)^{4cr/ \log (cr)} r_{i+1}\\
&\ll r_{i+1} P^\ve, 
\end{align*}
for $0 \leq i \leq n-1$.  We may conclude that
\begin{align*}
\prod_{i=\eta-1}^{n-1} q_{i}^{(i-\eta-1)/2}\ll
P^\ve 
\prod_{i=\eta-1}^{n-1} r_{i+1}^{(i-\eta-1)/2}\ll
P^\ve 
\prod_{j=\eta}^{n} r_{j}^{(j-\eta)/2},
\end{align*}
which therefore completes the argument in the case that $s_\infty\geq
0$. 

Suppose now that $s_\infty=-1$. Then 
\eqref{eq:induce} holds for integers $\eta$ in the shorter interval $[1,n]$. 
One easily checks in this case that 
$q_i \ll r_{i+1}P^\ve$ for $0\leq i \leq n-1$,
whence the upper bound is still satisfactory. 
It remains to deal with the case $s_{\infty}=-1$ and $\eta=0$. 
But this is exactly
the content of the Proposition \ref{lem:T1}, and so completes the
proof of Proposition \ref{lem:T2} subject to the resolution of
Proposition \ref{lem:T1}.

\section{Cubic exponential sums: the inductive base}\label{s:cubic_ib}

In this section we establish \prop \ref{lem:T1}. 
The essential ingredient in our estimation of $\TT(\al)$ will be 
an application of the Poisson
summation formula, in doing which we will draw inspiration from 
the second author's treatment of
cubic exponential sums in \cite{hb-10}.
Before embarking on the proof, we remind the reader of our convention 
concerning the value of $\ve$. 
Thus $\ve$ is a small positive parameter that is allowed to take 
different values at different parts of
the argument, and all of the implied constants are allowed to depend 
on $\ve$ without
further comment.  
Similarly, we will allow an implicit dependence on
$n$, on the constant $A$ that appears in the statement of \prop 
\ref{lem:T1}, and on the constants $c_n,c_{n,j}$ that feature in 
the definition of the set of weight function $\WW_n$.

Write $q=bc^2d$, where $b,d$ are given by \eqref{eq:bcd}. 
Our first step involves introducing complete
exponential sums modulo $q$. This will be achieved via an application
of Poisson summation.

\begin{lem}\label{lem:p_sum}
We have
$$
  \TT(a/q+z)=q^{-n}\sum_{\v\in\Z^n}T(a,q;\v)I(z;q^{-1}\v),
$$
where
\begin{equation}
  \label{eq:Sq}
T(a,q;\v):=\sum_{\y\bmod{q}}e_q(ag(\y)+\v.\y),
\end{equation}
and 
\begin{equation}
  \label{eq:Iq}
  I(z;\mbeta):=\int w(\x/P)e(z g(\x)-\mbeta.\x)\d\x.
\end{equation}
\end{lem}

\begin{proof}
Write $\x=\y+q\z$, for $\y \bmod{q}$, and $\al=a/q+z$. Then we obtain
\begin{align*}
\TT(\al)=\sum_{\y\bmod{q}}e_q(ag(\y))
\sum_{\z\in\Z^n}w((\y+q\z)/P)e(z g(\y+q\z)).
\end{align*}
An application of Poisson summation now yields
\begin{align*}
\TT(\al)
&=\sum_{\y\bmod{q}}e_q(ag(\y))\sum_{\v\in\Z^n}\int w((\y+q\z)/P)e(z
g(\y+q\z)-\v.\z)\d\z\\
&=q^{-n}\sum_{\v\in\Z^n}\sum_{\y\bmod{q}}e_q(ag(\y)+\v.\y)
\int w(\x/P)e(z g(\x)-\v.\x/q)\d\x.
\end{align*}
This completes the proof of the lemma.
\end{proof}

Lemma \ref{lem:p_sum} allows us to focus attention on a certain family
of complete exponential sums \eqref{eq:Sq} and integrals
\eqref{eq:Iq}. We begin with a treatment of the latter, when 
$\mbeta=q^{-1}\v$.

\begin{lem}\label{lem:Iq'}
Let $\ve>0$ and let $N\geq 1$. Then we have 
$$
I(z;q^{-1}\v) 
\ll_N |\v|^{- N}, 
$$
if $|\v|> HP^3$. Alternatively, when 
$|\v|\leq HP^3$, we have 
$$
I(z;q^{-1}\v) 
\ll_N 
P^{- N}+ \meas\big\{\x\ll P: ~|\v_0(\x) - \v|\leq
 P^\ve V\big\}, 
$$
where
\begin{equation}
  \label{eq:interval}
\v_0=\v_0(\x):=qz\nabla g(\x)
\end{equation}
and $V$ is given by \eqref{eq:V}.
\end{lem}

\begin{proof}
Define $g_P(\x):=P^{-3}g(P\x)$, for any $P\geq 1$. It is easily seen that
$$
I(z;q^{-1}\v)=P^n\int w(\x)e(z P^3 g_P(\x)-P\v.\x/q)\d\x=
P^nI_w(g_P,z P^3,P\v/q),
$$
in the notation of \eqref{eq:Iw}.  Moreover,  
$\|g_P\|=\|g\|_P\leq H$ by assumption. Hence Lemma~\ref{lem:Iq} implies that
$I(z;q^{-1}\v)$ is 
$$
\ll_N \frac{P^n}{R^N}+ P^n\meas\big\{\x\ll 1: ~\big|z P^2 
  \nabla g_P(\x) - \frac{\v}{q}\big|\leq \frac{R}{P} 
\max\{1,\sqrt{|z|HP^3}\}\big\},
$$
for any $R\gg 1$. We can take any fixed positive integer value for $N$
in the above, and we will change its value a number of times in what
follows, without further comment.  Suppose that $\v$ is contained in an
annulus 
$$
M<|\v|\leq 2M,
$$ 
for some $M>0.$
Then $|z P^2 \nabla g_P(\x) - \v/q|\gg M/q$ for any $\x\ll 1$, 
provided that $M\gg q|z|HP^2$.  On taking $R=M^{1/2}$ in our
estimate for $I(z;q^{-1}\v)$, we therefore deduce that
$
I(z;q^{-1}\v) \ll_N M^{-N}P^n,
$
for any $\v$ in the range $M<|\v|\leq 2M$, with 
$$
M\gg q|z|HP^2+
\frac{q M^{1/2}}{P} \max\{1,\sqrt{|z|HP^3}\}.
$$
Such an inequality clearly holds when $M\ge HP^3$,
by \eqref{eq:aqP} and \eqref{eq:z}.
It follows that $I(z;q^{-1}\v)\ll_N |\v|^{-N}$ for 
vectors $\v\in\Z^n$ with $|\v|> HP^3$.

Turning to the contribution from vectors $|\v|\leq HP^3$, we take
$R=P^\ve$ in our estimate for $I(z;q^{-1}\v)$. This implies that
$I(z;q^{-1}\v)$ is
$$
\ll_N \frac{1}{P^{\ve N}}+ \meas\big\{\x\ll P: 
~|z \nabla g(\x) - \frac{\v}{q}|\leq \frac{P^\ve}{P}
\max\{1,\sqrt{|z|HP^3}\}\big\}.
$$
On recalling the definitions \eqref{eq:V}, \eqref{eq:interval} of $V$ and
$\v_0(\x)$, this therefore suffices to complete the
proof of the lemma.
\end{proof}

Note that $\log H \ll \log P$. We may therefore combine Lemma
\ref{lem:Iq'} with Lemma \ref{lem:p_sum}, in
order to deduce that
\begin{align}  \label{eq:T2}
\TT(a/q+z)&\ll_N P^{-N}+q^{-n}\int_{\x\ll P}\sum_{|\v-\v_0|\leq P^\ve 
  V} 
|T(a,q;\v)| \d\x\nonumber\\
&\ll_N  P^{-N}+q^{-n}P^n\max_{\v_0}\sum_{|\v-\v_0|\leq P^\ve V}
|T(a,q;\v)|
\end{align}
for any $N\geq 1$. Here we have used the fact that if $N\ge 2n$, 
say, then  
$$
q^{-n}\sum_{|\v|\ge P^3}|\v|^{-N}q^n\ll P^{-N},
$$
for example.

At this point we should explain a key difference between our current
approach and that used in the second author's work \cite{hb-10}.  If
we were to follow this approach exactly, we would instead be led to
use the bound
$$
\TT(\al) \ll_N \frac{1}{P^{N}}+\frac{1}{q^{n}}\sum_{\v\ll
HP^{1+\ve}}|T(a,q;\v)|
\meas\{\x\ll P: ~|\v_0(\x)-\v|\leq P^{\ve}V\}.
$$
Estimating the measure above requires information about the size of
the Hessian $\det\ma{H}_{g_0}(\x)$.  Since we have not been able to get
appropriate estimates with suitable uniformity in $\|g_0\|$ we have
adopted the alternative procedure described above.  The new difficulty we
face is that our sum over $\v$ now runs over a small box
$|\v-\v_0|\leq P^\ve V$, rather than the larger one given by 
$\v\ll HP^{1+\ve}$.

It remains to study the average order of $T(a,q;\v)$,
as $\v$ ranges over a box with sides of length $V$,
centred upon a point $\v_0$.
Our investigation of this topic will draw inspiration from the
contents of  \cite[\S 6]{hb-10}, although a
number of key differences will become apparent.
First we need to establish multiplicativity in $q$ 
for the cubic exponential
sum $T(a,q;\v)$.

\begin{lem}\label{lem:T=mult}
Let $q=rs$, for coprime $r,s$. Let $\bar{r}, \bar{s}$ be integers such
that $r\bar{r}+s\bar{s}=1.$  Then we have
$$
T(a,rs;\v)=T(a\bar{s},r;\bar{s}\v)T(a\bar{r},s;\bar{r}\v).
$$
\end{lem}

\begin{proof}
This is standard, and so we will be brief.  As $\ma{r}$ ranges over
vectors modulo  $r$, and $\ma{s}$ ranges over such vectors 
modulo $s$, so $\y=r\bar{r}\ma{s}+s\bar{s}\ma{r}$ 
ranges over a
complete set of residues modulo $q=rs$.  Clearly 
$$
ag(\y)+\v.\y \equiv
r\bar{r}(ag(\ma{s})+\ma{s}.\v)+
s\bar{s}(ag(\ma{r})+\ma{r}.\v) \bmod{q},
$$
since $(r\bar{r})^j\equiv r\bar{r} \bmod{q}$ and 
$(s\bar{s})^j\equiv s\bar{s} \bmod{q}$, for any $j \geq 1$.
Hence it follows that
$$
e_q(ag(\y)+\v.\y)= 
e_s(\bar{r}(ag(\ma{s})+\ma{s}.\v))
e_r(\bar{s}(ag(\ma{r})+\ma{r}.\v)),
$$
which gives us the statement of Lemma \ref{lem:T=mult}. 
\end{proof}

For each $0\leq i\leq n$,  set
\begin{equation}
  \label{eq:def_bdi}
b_i:=\prod_{\colt{p^e\| b}{s_p=i-1}}p^e,\quad
d_i:=\prod_{\colt{p\mid d}{s_p=i-1}}p,
\end{equation}
where $b,d$ are given by \eqref{eq:bcd}.
Let $\bar{b},q^* \in \Z$ be such that $b\bar{b}+c^2dq^*=1$.
Then it follows from Lemmas \ref{deligne} and \ref{lem:T=mult} that
\begin{align}
  \label{eq:avS1}
T(a,q;\v)&=T(aq^*,b;q^*\v)T(a\bar{b},c^2d;\bar{b}\v)\nonumber\\
&\ll A^{\omega(b)} 
b^{n/2}(b_1b_2^2\cdots b_n^n)^{1/2}|T(a\bar{b},c^2d;\bar{b}\v)|
\nonumber\\
&\ll b^{n/2+\ve}(b_1b_2^2\cdots
b_n^n)^{1/2}|T(a\bar{b},c^2d;\bar{b}\v)|, 
\end{align}
for any fixed $\ve>0$.
We must now consider the size of the sum 
$T(a\bar{b},c^2d;\bar{b}\v)$, for given $a \bmod{c^2d}$ such that
$\hcf(a,c^2d)=1$, and given $\bar{b}\in\Z$ such that 
$b\bar{b}\equiv 1 \bmod{c^2d}$.  
We may assume henceforth that $\bar{b}$ is a positive integer, with
$1\leq \bar{b}<c^2d$.
 Given a vector $\x\in\Z^n$ and a positive integer $m$, let 
\begin{equation}
  \label{eq:M_d}
M_m(\x):=
\#\big\{\y \bmod{m}: \nabla^2 g(\x)\y \equiv \ma{0} \mod{m}\big\},
\end{equation}
and
\begin{equation}
  \label{eq:N_d}
N_m(\x):=
\#\big\{\y \bmod{m}: \ma{H}_{g_0}(\x)\y\equiv \ma{0}\mod{m}\big\},
\end{equation}
where $\ma{H}_{g_0}$ is given by \eqref{eq:hessian}.  Note that
$N_m(\x)$ is alternatively the number of $\y\bmod{m}$ for which
$\ma{H}_{g_0}(\y)\x\equiv \ma{0}\bmod{m}$, by \eqref{2.1a}.
If we write $g=g_0+f_2+f_1+f_0$, with each $f_i$ a
form of degree $i$, then it is clear that $\nabla^2 g =
\ma{H}_{g_0}+\ma{H}_{f_2}$. 
Our first task is to establish the
following result.

\begin{lem}
  \label{eq:avS2}
We have
$$
|T(a\bar{b},c^2d;\bar{b}\v)|\leq 
(c^2d)^{n/2}
\sum_{\colt{\ma{a}\bmod{c}}{c \mid (a\nabla
  g(\ma{a})+\v)}}M_{d}(\ma{a})^{1/2},
$$
where $M_d(\ma{a})$ is given by \eqref{eq:M_d}.
\end{lem}

\begin{proof}
Writing $\y = \ma{s}+cd \ma{t}$ in \eqref{eq:Sq}, we see that
\begin{align*}
T(a\bar{b},c^2d;\bar{b}\v)
&=\sum_{\ma{s}\bmod{cd}}
e_{c^2d}\big(a\bar{b}g(\ma{s})+\bar{b}\v.\ma{s}\big)
\sum_{\t\bmod{c}} e_{c}\big(\t.(a\bar{b}\nabla g(\ma{s})+\bar{b}\v)\big)\\ 
&=c^n
\sum_{\colt{\ma{s}\bmod{cd}}{
c\mid (a\bar{b}\nabla g(\ma{s})+\bar{b}\v)}}
e_{c^2d}\big(a\bar{b}g(\ma{s})+\bar{b}\v.\ma{s}\big). 
\end{align*}
Now write $\ma{s}=\ma{a}+c\ma{b}$, and note that the condition
on $\ma{s}$ in this sum implies that $a\nabla g(\ma{a})+\v =
c\,\ma{c}$, for some $\ma{c}\in \Z^n$.  Since $d\mid c$ it follows that
\begin{align*}
a\bar{b}g(\ma{s})+\bar{b}\v.\ma{s} 
&\equiv a\bar{b}g(\ma{a})+\bar{b}\v.\ma{a}
+\bar{b}c^2\big(\ma{b}.\ma{c} 
+\frac{a}{2}\ma{b}^T \nabla^2 g(\ma{a})\ma{b}\big) 
\mod{c^2d},
\end{align*}
whence 
\[
|T(a\bar{b},c^2d;\bar{b}\v)|=c^n\Big|
\sum_{\colt{\ma{a}\bmod{c}}{c \mid (a\nabla g(\ma{a})+\v)}}
S(\ma{a},\ma{c})\Big|,
\]
where
\[
S(\ma{a},\ma{c}):=
\sum_{\ma{b}\bmod{d}}e_{c^2d}\Big(
a\bar{b}g(\ma{a})+\bar{b}\v.\ma{a}+\bar{b}c^2\big(\ma{b}.\ma{c} 
+\frac{a}{2}\ma{b}^T \nabla^2 g(\ma{a})\ma{b} 
\big)\Big).\]
Moreover we have
\[|S(\ma{a},\ma{c})|\le \max_{\ma{c}\bmod{d}} |S_{\ma{a},\ma{c}}|,\]
with
$$
S_{\ma{a},\ma{c}}:=\sum_{\ma{b} \bmod{d}}
e_{d}\Big(\bar{b}\ma{b}.\ma{c}+ 
\frac{a\bar{b}}{2}\ma{b}^T \nabla^2 g(\ma{a})\ma{b}\Big).
$$
We estimate $S_{\ma{a},\ma{c}}$ by writing 
$$
|S_{\ma{a},\ma{c}}|^2=\sum_{\ma{b}_1,\ma{b}_2 \bmod{d}}
e_{d}\Big(\bar{b}\ma{c}.(\ma{b}_1-\ma{b}_2)+ 
\frac{a\bar{b}}{2}\big(
\ma{b}_1^T \nabla^2 g(\ma{a})\ma{b}_1-
\ma{b}_2^T \nabla^2 g(\ma{a})\ma{b}_2\big)
\Big).
$$
We write $\ma{b}_1=\ma{b}_2+\ma{b}_3$, and 
observe that
\begin{align*}
(\ma{b_2}+\ma{b}_3)^T\nabla^2 g(\ma{a})
&(\ma{b_2}+\ma{b}_3) -\ma{b}_2^T\nabla^2 g(\ma{a})\ma{b}_2\\
&=
\ma{b_2}.(2\nabla^2 g(\ma{a})\ma{b}_3) + \ma{b}_3^T\nabla^2 g(\ma{a})\ma{b}_3.
\end{align*}
We therefore obtain
$$
|S_{\ma{a},\ma{c}}|^2
\leq \sum_{\ma{b}_3 \bmod{d}}
\Big|
\sum_{\ma{b}_2\bmod{d}}e_{d}\big(
a\bar{b}\ma{b}_2.\nabla^2 g(\ma{a})\ma{b}_3\big)
\Big|
= d^n M_{d}(\ma{a}),
$$
in the notation of \eqref{eq:M_d}. 
It follows that 
$$
|T(a,c^2d;\bar{b}\v)| \leq 
(c^2d)^{n/2}\sum_{\colt{\ma{a}\bmod{c}}{c \mid (a\nabla
  g(\ma{a})+\v)}}M_{d}(\ma{a})^{1/2}\\
$$
which completes the proof of Lemma \ref{eq:avS2}.
\end{proof}

Recall the assumptions \eqref{eq:aqP} and \eqref{eq:z} on $a,q,z$. 
We can now combine \eqref{eq:T2} and \eqref{eq:avS1} with Lemma
\ref{eq:avS2}, and obtain the following conclusion.

\begin{lem}\label{combo}
Define
\begin{equation}
  \label{eq:SV}
  \mcal{S}(V,a)=\mcal{S}(V,a;\v_0,c,d):=
\sum_{|\ma{v}-\ma{v}_0|\leq V}\sum_{\colt{\ma{a}\bmod{c}}{c \mid (a\nabla
  g(\ma{a})+\v)}}M_{d}(\ma{a})^{1/2}.
\end{equation}
Then we have
$$
\TT(a/q+z)\ll_N P^{-N}+q^{-n/2}(b_1b_2^2\cdots b_n^n)^{1/2} 
P^{n+\ve}\max_{\v_0}\mcal{S}(P^{\ve}V,a).
$$
\end{lem}

We must now make a closer examination of the sum $\mcal{S}(V,a)$.
One of the ingredients that goes into this investigation is the 
average order of the function
$N_{m}(\ma{r})^{1/2}$, where $N_{m}(\ma{r})$ is given by
\eqref{eq:N_d}. Specifically we
will need the following result.

\begin{lem} \label{eq:avS5}
Let $R\geq 1$ and let $m\in \N$. Then we have 
$$
\sum_{{|\ma{r}|\leq R}} N_{m}(\ma{r})^{1/2}
\ll
m^{n/2}\min \Big\{R, 
~ \Big(1+ 
\frac{HR^3}{m}\Big)^{1/2}\Big\}^{n}.
$$
\end{lem}
\begin{proof}
To start with, it is trivial to see that
$$
\sum_{{|\ma{r}|\leq R}} N_{m}(\ma{r})^{1/2}
\ll
 m^{n/2}R^n,
$$
which is satisfactory for the lemma.

To obtain an alternative estimate, we suppose first that 
$RH \leq m$. 
We think of $\y$ in the definition \eqref{eq:N_d} as
running over $(0,m]^n$, 
and split this region into $K^n$ 
subcubes of side $m/K$, where $K$ is a positive integer parameter at our
disposal.  If $\y_1$ and $\y_2$ are both solutions to
$\ma{H}_{g_0}(\ma{r})\y\equiv \ma{0}\bmod{m}$, lying in the same subcube,
then $\y_3=\y_2-\y_1$ is also a solution, and lies in $(-m/K,m/K)^n$.
We conclude that
$$
N_m(\ma{r})\le K^n\#\{\y\in(-m/K,m/K)^n:
\ma{H}_{g_0}(\ma{r})\y\equiv \ma{0}\mod{m}\}.
$$
We now choose $K$ of order $RH$ so that $\y\in(-m/K,m/K)^n$
implies that $|(\ma{H}_{g_0}(\ma{r})\y)_i|<m$ for each $i$. 
Here we have used \eqref{eq:cots} to deduce that the matrix
$\ma{H}_{g_0}(\ma{r})$ has entries of order $O(RH)$.

It now follows that 
$$
N_m(\ma{r})\ll (RH)^n\#\{\y\in(-m/K,m/K)^n:\ma{H}_{g_0}(\ma{r})\y=\ma{0}\}.
$$
However the condition $\ma{H}_{g_0}(\ma{r})\y=\ma{0}$ restricts $\y$
to a linear space of dimension $n-\rho(\ma{r})$, where
$
\rho(\ma{r}):= \rank \ma{H}_{g_0}(\ma{r}).
$
Thus \eqref{eq:aff_triv} implies that
$$
N_m(\ma{r})\ll (RH)^n\Big(\frac{m}{RH}\Big)^{n-\rho(\ma{r})}, 
$$
since $RH\leq m.$  
Taken together, \eqref{eq:aff_triv} and Lemma \ref{lem:homrank} now show that
\begin{align*}
\sum_{{|\ma{r}|\leq R}} N_{m}(\ma{r})^{1/2}
&\ll(RH)^{n/2} 
\sum_{t=0}^n
\Big(\frac{m}{RH}\Big)^{(n-t)/2}\#\{|\ma{r}|\leq R: \rho(\ma{r})=t\}\\
&\ll(RH)^{n/2} 
\sum_{t=0}^n
\Big(\frac{m}{RH}\Big)^{(n-t)/2}R^{t}\\
&\ll(RH)^{n/2}\Big(R^n +
\Big(\frac{m}{RH}\Big)^{n/2}\Big)\\
&= m^{n/2}\Big(1+
\Big(
\frac{R^3H}{m}
\Big)^{n/2}\Big).
\end{align*}
Finally, if $RH >m$, then we trivially have
\begin{align*}
\sum_{{|\ma{r}|\leq R}} N_{m}(\ma{r})^{1/2}
\ll m^{n/2} R^n \ll m^{n/2}\Big(1+
\Big(\frac{R^3H}{m}
\Big)^{n/2}\Big).
\end{align*}
This completes the proof of Lemma \ref{eq:avS5}.
\end{proof}

We are now ready to proceed with our
analysis of $\mcal{S}(V,a)$, for which we will provide two
alternative estimates.  It will be convenient to set
\begin{equation}
  \label{eq:Di}
  D:=d_1d_2^2\cdots d_n^n,
\end{equation}
in what follows, where $d_1,\ldots,d_n$ are given by \eqref{eq:def_bdi}.
We begin by recording the following simple estimate.

\begin{lem}\label{lem:avS3}
We have
$$
  \sum_{\ma{a} \bmod{d}} M_{d}(\ma{a})
\ll d^{n+\ve} D.
$$
\end{lem}
\begin{proof}
The result is trivial if $d=1$. Suppose that $d>1$.
Since $d$ is square-free we may write $d=p_1\cdots p_r$, for distinct 
primes $p_1,\ldots, p_r$.
Hence we may combine Lemma \ref{lem:affhess} with \eqref{eq:n-tcor}, 
to deduce that
\begin{align*}
\sum_{\ma{a} \bmod{d}} M_{d}(\ma{a}) 
&=
\#\big\{\ma{a}, \ma{b} \bmod{d}: \nabla^2 g(\ma{a})\ma{b} 
\equiv \ma{0} \mod{d}\big\}\\
&\ll A^r \prod_{i=1}^r p_i^{n+1+s_{p_i}}
\ll d^{n+\ve} d_1d_2^2\cdots d_n^n,
\end{align*}
in the notation of \eqref{eq:def_bdi}. Here $A$ is a constant
depending only on $n$, and we have used 
the fact that $A^{\omega(d)}
\ll_{A} d^{\ve}$. This completes the proof of the lemma.
\end{proof}

We are now ready to record our first bound for 
$\mcal{S}(V,a)$, as given by \eqref{eq:SV}.

\begin{lem}\label{eq:s_1}
We have
$$
\SSS(V,a)\ll c^{\ve}
D^{1/2}V^n\Big(1+\frac{c}{V}\Big)^{n/2}
\Big(1+\min\Big\{\frac{c}{V},
\frac{c^{2/3}H^{1/3}}{V}\Big\}\Big)^{n/2},
$$
where $D$ is given by \eqref{eq:Di}.
\end{lem}

\begin{proof}
By Cauchy's inequality, we have $\SSS(V,a)\le
\SSS_1(V,a)^{1/2}\SSS_2(V,a)^{1/2}$, where
$$
\SSS_1(V,a):=
\sum_{|\ma{v}-\ma{v}_0|\leq V}\sum_{\colt{\ma{a}\bmod{c}}{c \mid (a\nabla
  g(\ma{a})+\v)}}M_{d}(\ma{a})
$$
and
$$
\SSS_2(V,a):=
\sum_{|\ma{v}-\ma{v}_0|\leq V}\sum_{\colt{\ma{a}\bmod{c}}{c \mid (a\nabla
  g(\ma{a})+\v)}}1.
$$
We begin by considering $\SSS_1(V,a)$.  We have
\begin{align}\label{eq:sK_large}
\SSS_1(V,a)\le&\sum_{\ma{a}\bmod{c}}M_{d}(\ma{a})\#\{\v:
|\ma{v}-\ma{v}_0|\leq V,\,c \mid a\nabla g(\ma{a})+\v\}\nonumber\\
\ll&\sum_{\ma{a}\bmod{c}}M_{d}(\ma{a})\Big(1+\frac{V}{c}\Big)^n\nonumber\\
=&\Big(\frac{c}{d}\Big)^n\Big(1+\frac{V}{c}\Big)^n
\sum_{\ma{a}\bmod{d}}M_{d}(\ma{a})\nonumber\\
\ll&\Big(\frac{c}{d}\Big)^n\Big(1+\frac{V}{c}\Big)^nd^{n+\ve}D\nonumber\\
=&D(V+c)^nd^{\ve},
\end{align}
by Lemma \ref{lem:avS3}. 

We turn now to $\SSS_2(V,a)$, for which we will show that
\begin{equation}
  \label{eq:sK_small}
\SSS_2(V,a) \ll c^{\ve} \big(V+\min\{c,(c^2H)^{1/3}\}\big)^{n}
\end{equation}
when $c\gg V$.  We begin by noting that
$$
\SSS_2(V,a)\ll
\sum_{\ma{a}\bmod{c}}\sum_{\colt{\v}{c \mid (a\nabla g(\ma{a})+\v)}} 
\exp(-\|\ma{v}-\ma{v}_0\|^2 V^{-2}), 
$$
where $\|\ma{z}\|:=\sqrt{z_1^2+\cdots+z_n^2}$ denotes the Euclidean 
norm on $\R^n$. 
On appealing to the Poisson summation formula we find that the inner
sum is
$$
\pi^{n/2}\Big(\frac{V}{c}\Big)^n\sum_{\ma{r}\in\Z^n}
e_c\big(\ma{r}.(\v_0+a\nabla g(\ma{a}))\big)\exp(-\pi^2\|\ma{r}\|^2V^2/c^2),
$$
whence
$$
\SSS_2(V,a)\ll\Big(\frac{V}{c}\Big)^n\sum_{\ma{r}\in\Z^n}
\exp(-\pi^2\|\ma{r}\|^2V^2/c^2)|\tau_{\ma{r}}|, 
$$
with
$$
\tau_\ma{r}:=\sum_{\ma{a}\bmod{c}} e_{c}(a\ma{r}.\nabla g (\ma{a})).
$$
Terms with $|\ma{r}|\gg c(\log c)/V$ trivially contribute
\begin{align*}
&\ll V^n\Big(\sum_{r\in\Z}\exp(-\pi^2r^2V^2/c^2)\Big)^{n-1}
\Big(\sum_{|r|\gg c(\log c)/V}\exp(-\pi^2r^2V^2/c^2)\Big)\\
&\ll V^n(c/V)^{n-1}\big(c\exp(-\log^2c)/V\big)\\ 
&\ll 1,
\end{align*}
whence
\begin{equation}\label{eq:intro_tau2}
\SSS_2(V,a)\ll 1+\Big(\frac{V}{c}\Big)^n\sum_{\ma{r}\ll c(\log c)/V}
|\tau_{\ma{r}}|.
\end{equation}
We now observe that
\begin{align*}
  |\tau_\ma{r}|^2
&=
\sum_{\ma{a}_1,\ma{a}_2\bmod{c}}
e_{c}\big(a\ma{r}.(\nabla g (\ma{a}_1)-\nabla g (\ma{a}_2))\big).
\end{align*}
If we write $\ma{a}_1=\ma{a}_2+ \ma{a}_3$, we find that
$\nabla g (\ma{a}_2+ \ma{a}_3)-\nabla g (\ma{a}_2)$
is equal to $\ma{H}_{g_0}(\ma{a}_2) \ma{a}_3$ plus
a term that is independent of $\ma{a}_2$. 
It therefore follows from \eqref{2.1a} that 
\begin{align*}
  |\tau_\ma{r}|^2
&\leq \sum_{\ma{a}_3\bmod{c}}\Big|
\sum_{\ma{a}_2\bmod{c}}
e_{c}\big(\ma{a}_2.(a\ma{H}_{g_0}(\ma{r}) \ma{a}_3)\big)\Big|\\
&\ll  c^nN_{c}(\ma{r}),
\end{align*}
in the notation of \eqref{eq:N_d}.
Substituting this into \eqref{eq:intro_tau2}, we therefore conclude
that
$$
\SSS_2(V,a)\ll 1+\frac{V^n}{c^{n/2}}
\sum_{\ma{r}\ll c(\log c)/V} N_c(\ma{r})^{1/2},
$$
An application of Lemma \ref{eq:avS5}, with $R\gg c(\log c)/V$ 
and $m=c$, now yields
\begin{align*}
\SSS_2(V,a) \ll& c^{\ve}V^{n}\min \Big\{\frac{c}{V}, 
~ \Big(1+ \frac{c^2H}{V^3}\Big)^{1/2}\Big\}^{n}\\
\ll&  c^{\ve} V^{n}\Big(1+\frac{c}{V}\min \Big\{1,~ 
\frac{H^{1/2}}{V^{1/2}}\Big\}\Big)^n\\
\ll& c^{\ve}  \Big(V+c\min \Big\{1,~ 
\frac{H^{1/2}}{V^{1/2}}\Big\}\Big)^{n},
\end{align*}
when $c\gg V$. At this point we make the observation that the
quantity $\SSS_2(V,a)$ can only be made larger by
increasing the size of $V$.  On writing 
$$
V_0:=V+\min\{c,(c^2H)^{1/3}\},
$$
we note in particular that $V\leq V_0 \ll c$. It therefore follows
from the above estimate that
\begin{align*}
\SSS_2(V,a) &\leq \SSS_2(V_0,a)\\
&\ll c^{\ve}\Big(V_0+c\min \Big\{1,~ 
\frac{H^{1/2}}{c^{1/2}}
+\frac{H^{1/3}}{c^{1/3}}\Big\}\Big)^{n}\\
&\ll c^{\ve} V_0^{n}.
\end{align*}
This completes the proof of \eqref{eq:sK_small}.

Let us continue to adopt the notation for $V_0$ introduced above.
We are now in a position to combine \eqref{eq:sK_large} and 
\eqref{eq:sK_small} to deduce that
\begin{align*}
\SSS(V,a)&\le\SSS_1(V,a)^{1/2}\SSS_2(V,a)^{1/2}\\
&\ll \big(d^{\ve} D(V+c)^n\big)^{1/2}
\big(c^{\ve} \big(V+\min\{c,(c^2H)^{1/3}\}\big)^{n}\big)^{1/2}.
\end{align*}
This suffices for Lemma \ref{eq:s_1} when $c\ge V$.

In the remaining case $c\le V$, we bound $\SSS_2(V,a)$ trivially as
\begin{align*}
\SSS_2(V,a)&\ll
\sum_{\ma{a}\bmod{c}}\#\{\v:\,|\v-\v_0|\le V,\,\v\equiv -a\nabla g(\ma{a})
\bmod{c}\}\\
&\ll c^n(1+V/c)^n\\
&\ll V^n.
\end{align*}
Since \eqref{eq:sK_large} yields $\SSS_1(V,a)\ll d^{\ve} DV^n$ we see 
that Lemma \ref{eq:s_1} follows for $c\le V$ too. 
\end{proof}

The following result provides an alternative estimate for $\SSS(V,a)$,
and follows from a rather straightforward
modification to the proof of Lemma~\ref{eq:s_1}.

\begin{lem}\label{eq:s_2}
We have
$$
\SSS(V,a)
\ll c^{\ve}D^{1/2}
V^{n}\Big(1+ \frac{Hc^2d}{V^3}\Big)^{n/2},
$$
where $D$ is given by \eqref{eq:Di}.
\end{lem}
\begin{proof}
When  $c<V$ this follows directly from Lemma \ref{eq:s_1}. When $c\geq V$ we
follow the proof of \eqref{eq:intro_tau2}, but apply the method to
$\SSS(V,a)$ directly.  This yields
\begin{equation}
  \label{eq:intro_tau'}
\SSS(V,a)\ll 1+\Big(\frac{V}{c}\Big)^n
\sum_{\ma{r}\ll c(\log c)/V}  |\sigma_{\ma{r}}|,
\end{equation}
where now
\begin{align*}
\sigma_\ma{r}&:=
\sum_{\ma{a}\bmod{c}} e_{c}(a\ma{r}.\nabla g (\ma{a}))
M_{d}(\ma{a})^{1/2}\\
&=
\sum_{\ma{b}\bmod{d}}
M_{d}(\ma{b})^{1/2}
\sum_{\colt{\ma{a}\bmod{c}}{\ma{a}\equiv \ma{b} \bmod{d}}}
e_{c}(a\ma{r}.\nabla g (\ma{a})).
\end{align*}
On combining Lemma \ref{lem:avS3} with an application of Cauchy's
inequality, we deduce that
\begin{align*}
  |\sigma_\ma{r}|^2
&\leq
\Big(\sum_{\ma{b}\bmod{d}}
M_{d}(\ma{b})\Big)\Big(
\sum_{\ma{b}\bmod{d}}\Big|
\sum_{\colt{\ma{a}\bmod{c}}{\ma{a}\equiv \ma{b} \bmod{d}}}
e_{c}(a\ma{r}.\nabla g (\ma{a}))\Big|^2\Big)\\
&\ll c^{\ve}d^{n}D\Big|
\sum_{\colt{\ma{a}_1,\ma{a}_2\bmod{c}}{\ma{a}_1\equiv \ma{a}_2 \bmod{d}}}
e_{c}\big(a\ma{r}.(\nabla g (\ma{a}_1)-\nabla g (\ma{a}_2))\big)\Big|.
\end{align*}
We now write $\ma{a}_1=\ma{a}_2+d \ma{a}_3$, and find that
$\nabla g (\ma{a}_2+d \ma{a}_3)-\nabla g (\ma{a}_2)$
is equal to $d\ma{H}_{g_0}(\ma{a}_3)\ma{a}_2$ plus terms that are 
independent of $\ma{a}_2$. 
It easily follows that 
$$
  |\sigma_\ma{r}|^2
\ll  c^{n+\ve}d^n D N_{c/d}(\ma{r}).
$$
Substituting this into \eqref{eq:intro_tau'}, and applying
Lemma \ref{eq:avS5} with $R\gg c(\log c)/V$ and $m=c/d$, 
we therefore conclude the proof of Lemma \ref{eq:s_2}.
\end{proof}

We proceed by noting that
$$
b_1b_2^2\cdots b_n^n D = b_1d_1b_2^2d_2^2\cdots b_n^nd_n^n =
r_1r_2^2\cdots r_n^n,
$$
where $r_1,\ldots, r_n$ are given by \eqref{eq:qi}.
Putting together
Lemmas \ref{combo}, \ref{eq:s_1} and \ref{eq:s_2}, it therefore 
follows that 
\begin{align*}
\TT(a/q+z)&\ll q^{-n/2+\ve}(r_1r_2^2\cdots r_n^n)^{1/2}P^nV^n\\
&\times
\min\Big\{
\Big(1+\frac{c}{V}\Big)\Big(1+\min\Big\{\frac{c}{V},
\frac{c^{2/3}H^{1/3}}{V}\Big\}\Big),
1+  \frac{c^2dH}{V^3}\Big\}^{n/2}.
\end{align*}
Now it is easy to see that 
$$
\Big(1+\frac{c}{V}\Big)\Big(1+\min\Big\{\frac{c}{V},
\frac{c^{2/3}H^{1/3}}{V}\Big\}\Big) \ll
1+\frac{c}{V}+\min\Big\{\frac{c^2}{V^2},
\frac{c^{5/3}H^{1/3}}{V^2}\Big\},
$$
whence
\begin{align*}
\min\Big\{
\Big(1+\frac{c}{V}\Big)&\Big(1+
\min\Big\{\frac{c}{V},
\frac{c^{2/3}H^{1/3}}{V}\Big\}\Big),
1+  \frac{c^2dH}{V^3}\Big\}\\
&\ll 1+ 
\min\Big\{\frac{c}{V},\frac{c^2dH}{V^3} \Big\}+
\min\Big\{
\frac{c^2}{V^2},
\frac{c^{5/3}H^{1/3}}{V^2}, \frac{c^2dH}{V^3}
\Big\}.
\end{align*}
We have therefore established that
\begin{equation}\label{eq:T33}
\TT(a/q+z)\ll q^{-n/2}(r_1r_2^2\cdots r_n^n)^{1/2}P^{n+\ve}V^n
\big(1+M_1(V)+M_2(V)\big)^{n/2}, 
\end{equation}
where
$$
M_1(V):=\min\Big\{\frac{c}{V},\frac{c^2dH}{V^3} \Big\},\quad
M_2(V):= \min\Big\{
\frac{c^2}{V^2},
\frac{c^{5/3}H^{1/3}}{V^2}, \frac{c^2dH}{V^3}
\Big\}.
$$

We are now ready to complete our proof of Proposition~\ref{lem:T1}.
Suppose first that $c<V$, where 
$V$ is given by \eqref{eq:V}.
Then it follows from \eqref{eq:T33} that  
\begin{align*}
 \TT(a/q+z)
&\ll q^{-n/2}r_1^{1/2}\cdots r_n^{n/2}P^{n+\ve} V^n.
\end{align*}
This is satisfactory for Proposition \ref{lem:T1}.
Suppose now that $c\geq V$, and set 
$$ 
V_1:=V+(c^2d H)^{1/3}.
$$
In particular we have $V\leq V_1$. 
It should be clear from Lemma \ref{combo} that the upper bound in
\eqref{eq:T33} remains valid when $V$ is replaced by anything that
exceeds it. Hence 
\begin{equation}\label{trefoil}
\TT(a/q+z)\ll q^{-n/2}(r_1r_2^2\cdots r_n^n)^{1/2}P^{n+\ve}V_1^n, 
\end{equation}
since $M_i(V_1)\leq 1$ for $i=1,2$. 
We will obtain an alternative estimate for $\TT(a/q+z)$ by taking
$$
1+M_1(V)+M_2(V)\leq 1+
\frac{c}{V}+\frac{c^{5/3}H^{1/3}}{V^2}
$$
in \eqref{eq:T33}.
Still under the assumption that $c\geq V$, we deduce that 
\begin{align*}
\TT(a/q+z)&\ll q^{-n/2} r_1^{1/2}\cdots r_n^{n/2}P^{n+\ve}  
\big(V+c^{1/2}V^{1/2}+c^{5/6}H^{1/6}\big)^{n}.
\end{align*}
Taken together with \eqref{trefoil}, we find that
\begin{align*}
\TT(a/q+z)
\ll q^{-n/2} r_1^{1/2}\cdots r_n^{n/2} P^{n+\ve}W^{n}, 
\end{align*}
where $W$ is given by \eqref{eq:W}. 
This completes the proof of Proposition \ref{lem:T1}.

\section{Estimating $S(\al)$: Weyl differencing}\label{s:weyl}

In this section and the next, our aim is to estimate 
the quartic exponential sum \eqref{eq:S}. The results obtained will
form a key ingredient in our application of the circle method to
estimate \eqref{eq:start}, particularly in the context of the minor 
arcs.  We will need to say a few words about
the function $\omega:\R^n\rightarrow \R_{\geq 0}$ that appears in the
definition of $S(\al)$. In \S \ref{s:vdc} we will need to suppose that
$\omega\in \WW_n$, as defined before Lemma \ref{lem:Iq}. 
In the present section, 
which is dedicated to describing the route taken by Birch,
there is no need to be so restrictive. Thus we will suppose only that
$\omega\in\WW_n\cup\{\chi\}$, where $\chi:\R^n\rightarrow \R_{\geq 0}$ 
denotes the
characteristic function on $(0,1]^n$.
Throughout this section, it will be convenient to set
\begin{equation}
  \label{eq:sigma}
  \sigma:=\dim \sing_\Q(X).
\end{equation}
Thus $\sigma$ is an integer in the interval $[-1,n-3]$, with the usual
convention that $\sigma=-1$ if and only if $X$ is non-singular.

The central idea in Birch's approach involves Weyl differencing.
The first step in this process produces the bound
\begin{equation}\label{weyl1}
|S(\alpha)|^2\ll \sum_{\w\ll P}\Big|\sum_{\x\in\Z^n} 
\omega\big((\x+\w)/P\big)\omega(\x/P)
e\big(\al (F(\x+\w)-F(\x))\big)\Big|.
\end{equation}
A further application of Cauchy's inequality now yields
$$
|S(\alpha)|^4\ll P^n\sum_{\w,\x\ll P}\Big|\sum_{\y\in \Z^n}
\omega_{\w,\x}(\y)
e\big(\alpha F(\w,\x;\y)\big)\Big|,
$$
where 
$$
F(\w,\x;\y):=F(\w+\x+\y)-F(\w+\y)-F(\x+\y)+F(\y)
$$
and 
$$
\omega_{\w,\x}(\y)=\omega\big((\w+\x+\y)/P\big)
\omega\big((\w+\y)/P\big)\omega\big((\x+\y)/P\big)\omega(\y/P).
$$
We now repeat this procedure, obtaining
\begin{equation}\label{22-weyl3}
|S(\alpha)|^8\ll P^{4n}\sum_{\w,\x,\y\ll P}\Big|\sum_{\z\in\Z^n}
\omega_{\w,\x,\y}(\z)
e\big(\alpha F(\w,\x,\y;\z)\big)\Big|,
\end{equation}
where now
\begin{align*}
F(\w,\x,\y;\z):=&F(\w+\x+\y+\z)-F(\w+\x+\z)-F(\w+\y+\z)\\&
\quad -F(\x+\y+\z)+F(\w+\z)+F(\x+\z)+F(\y+\z)\\
& \quad -F(\z),
\end{align*}
and $\omega_{\w,\x,\y}(\z)$ is defined in the obvious way.

We now recall the definition \eqref{eq:tri} of the trilinear forms
$L_i(\w;\x;\y)$, for $1\le i\le n$. It is not hard to see that
$$
F(\w,\x,\y;\z)=\sum_{i=1}^n z_iL_i(\w;\x;\y)+\Phi(\w,\x,\y),
$$
where $\Phi(\w,\x,\y)$ is independent of $\z$.  
It therefore follows from \eqref{22-weyl3} that
\begin{align*}
|S(\alpha)|^8
&\ll 
P^{4n}\sum_{\w,\x,\y\ll P}\Big|\sum_{\z\in \Z^n}
\omega_{\w,\x,\y}(\z)e\Big(\alpha \sum_{i=1}^n z_iL_i(\w;\x;\y)\Big)\Big|.
\end{align*}
If $\omega=\chi$, the characteristic function on $(0,1]^n$, 
then we have
\begin{equation}\label{22-weyl4}
|S(\alpha)|^8 
\ll P^{4n} \sum_{\w,\x,\y\ll P}\,
\prod_{i=1}^n
\min\{P, \| \alpha L_i(\w;\x;\y)\|^{-1}\}.
\end{equation}
If instead $\omega\in \WW_n$, then an application of partial
summation yields the same inequality. This estimate 
corresponds to \cite[Lemma 2.1]{birch} in the case
$R=1$ and $d=4$, with $\Phi_J(\al;\x^{(1)},\x^{(2)},\x^{(3)})=
\al L_J(\x^{(1)};\x^{(2)};\x^{(3)})$.

We proceed to define the quantity
$$
N(\alpha,P):=\#\big\{\w,\x,\y \ll P:~
\|\alpha L_i(\w;\x;\y)\| <P^{-1}~ \forall i\leq n\big\}.
$$
It is then a simple matter to deduce that
$$
\sum_{\w,\x,\y\ll P}\,
\prod_{i=1}^n
\min\{P, \| \alpha L_i(\w;\x;\y)\|^{-1}\}
\ll (P\log P)^n N(\alpha,P),
$$
as in the proof of \cite[Lemma 13.2]{dav-book}.
On inserting this into \eqref{22-weyl4}, it therefore follows that
\begin{equation}\label{22-smee}
|S(\alpha)|^8  \ll 
P^{5n}(\log P)^n N(\alpha,P),
\end{equation}
which corresponds to \cite[Lemma 2.2]{birch}.  In particular it is
clear from the trivial upper bound $N(\alpha, P)\ll P^{3n}$
that we have lost very little in formulating \eqref{22-smee}.
In order to handle the quantity $N(\alpha,P)$ we will employ the following
result, which is due to Davenport \cite[Lemma 12.6]{dav-book}.

\begin{lem}\label{22-davlem}
Let $\ma{L}$ be a real symmetric $n\times n$ matrix.  Let $A,c>0$
be real, and let
$$
N(Z) := \#\{\ma{u}\in\Z^n: ~|\ma{u}|\leq c AZ,~ \|(\ma{L}\ma{u})_i\| 
<A^{-1}Z~ \forall i\le n\}.
$$
Then, if $0<Z_1\le Z_2\le 1$, we have
$$
N(Z_2)\ll_c\Big(\frac{Z_2}{Z_1}\Big)^nN(Z_1).
$$
\end{lem}

The version of Lemma \ref{22-davlem} established by Davenport
corresponds to taking $c=1$. An inspection of the proof reveals that
the only difference involved in taking $c>0$ to be arbitrary is that
the implied constant in the upper bound for $N(Z_2)$ is allowed to
depend on $c$.

The idea is now to apply this result three times, in order to reduce 
the analysis of $N(\alpha,P)$ to a problem involving the system of
equations $L_i(\w;\x;\y)=0$, 
for $1\leq i\leq n$.   To begin with one takes the matrix $\ma{L}$ in
Lemma \ref{22-davlem} to be given by
$(\ma{L}\y)_i=L_i(\w;\x;\y)$.  Choosing $A=P$,
we have
$$
N(1)=\#\{\y\in\Z^{n}: |\y|\leq cP,~
\|\alpha L_i(\w;\x;\y)\| <P^{-1} ~\forall i\le n\},
$$
and the lemma implies that
$N(1)\ll_c Z^{-n}N(Z)$ for any $0<Z\le 1$ and any $c>0$.  It follows that
there exists a positive absolute constant $c=O(1)$ such that 
$$
N(\alpha,P)
\ll
Z^{-n}\#\Big\{(\w,\x,\y)\in\Z^{3n}:
\begin{array}{l}
|\w|,|\x|\leq cP, ~|\y|\leq cZP,\\
\|\alpha L_i(\w;\x;\y)\|<ZP^{-1}~ \forall i\le n
\end{array}
\Big\}.
$$
Rather than using this estimate directly, Birch permutes the 
r\^oles of $\w,\x$ and $\y$ in
the above.  Taking $\ma{L}$ to be the matrix given by
$(\ma{L}\x)_i=L_i(\w;\x;\y)$, one applies Lemma \ref{22-davlem} with 
$A=PZ^{-1/2}$, $Z_2=Z^{1/2}$ and $Z_1=Z^{3/2}$, in order to deduce that
$$
N(\alpha,P)
\ll
Z^{-2n}\#\Big\{(\w,\x,\y)\in\Z^{3n}:
\begin{array}{l}
|\w|\leq cP, ~|\x|,|\y|\leq cZP,\\
\|\alpha L_i(\w;\x;\y)\|<Z^2P^{-1}~ \forall i\le n
\end{array}
\Big\}.
$$
Finally, it remains to use Lemma \ref{22-davlem} to shrink the size of
the box that $\w$ lies in.
Taking $\ma{L}$ to be the matrix given by
$(\ma{L}\w)_i=L_i(\w;\x;\y)$, therefore, we apply Lemma \ref{22-davlem} with 
$A=PZ^{-1}$, $Z_2=Z$ and $Z_1=Z^{2}$, to conclude that
$$
N(\alpha,P)
\ll
Z^{-3n}\#S(cZP,Z^{-3}P).
$$
Here we have set 
$$
S(R,Q):=\Big\{(\w,\x,\y)\in\Z^{3n}:
\begin{array}{l}
|\w|,|\x|,|\y|\leq R,\\
\|\alpha L_i(\w;\x;\y)\|<Q^{-1}~ \forall i\le n
\end{array}
\Big\},
$$
for any $R,Q>0$.

The idea now is to find conditions on $\al$ and $Z$ under which 
\begin{equation}\label{23-bileq}
L_i(\w;\x;\y)=0, \quad (1\leq i\leq n),
\end{equation}
for every $(\w,\x,\y)\in S(cZP,Z^{-3}P)$.  
Rather than pursuing Birch's formulation of this
particular step, which requires the introduction of a slightly sparser
set of major arcs, we have decided to take advantage
of the second author's 
recent contribution to the topic \cite{14}.  This
includes a proof of the following simple result \cite[Lemma 2.3]{14}.

\begin{lem}\label{23-approx}
Let $M>0$ and let $\alpha=a/q+z$, with $|z|\leq
(2qM)^{-1}$. Suppose that $m\in \Z$ is such that $|m|\leq M$ and 
$||\alpha m||<Q^{-1}$ for some $Q\geq 2q$.  Then $q|m$.  In
particular we will have $m=0$ if in addition we have either 
$M<q$ or $|z| > (qQ)^{-1}$.
\end{lem}

We will use Lemma \ref{23-approx} to reduce
our consideration to the system of trilinear equations \eqref{23-bileq}.
Write $f=24 \sum|f_{ijk\ell}|$, where $f_{ijk\ell}$  are the coefficients
of $F$ in \eqref{eq:F}, and suppose that $\al=a/q+z$.  Then on
choosing $Z$ to satisfy the conditions
\begin{equation}\lab{23-Z1}
0<Z\leq 1,\quad  
Z^3\leq (2c^3fq|z|P^3)^{-1},\quad 
Z^3\leq P/(2q),
\end{equation}
and
\begin{equation}\lab{23-Z2}
Z^3<\max\Big\{\frac{q}{c^3fP^3},q|z|P\Big\},
\end{equation}
we may make Lemma \ref{23-approx} applicable, and therefore deduce that
\begin{equation}\lab{4-N}
N(\alpha,P)\ll Z^{-3n}\#T(cZP),
\end{equation}
where $c=O(1)$ as usual, and
$$
T(R):=\{(\w,\x,\y)\in(\Z\cap[-R,R])^{3n}:~ L_i(\w;\x;\y)=0~\forall
i\le n\}
$$
for any $R>0$.

We are now led to study the density of integer solutions to the system
of equations \eqref{23-bileq}.  Arguing as in the proof of 
\eqref{eq:diag}, which is a special
case of  \cite[Lemma~3.3]{birch}, it is easy to see that the 
variety cut out by \eqref{23-bileq} in $\A^{3n}$ has
dimension at most $2n+\sigma+1$, where $\sigma$ is defined in
\eqref{eq:sigma}. 
An application of \eqref{eq:aff_triv} now reveals  that
$$
\#T(R)\ll R^{2n+\sigma+1},
$$
for any $R\geq 1$. This is clearly best possible when $F$ is
non-singular.  We may now insert this bound into \eqref{4-N} to conclude that
\begin{equation}\lab{4-mental}
N(\alpha,P)\ll Z^{-n+\sigma+1}P^{2n+\sigma+1},
\end{equation}
provided that $Z \geq P^{-1}$.  This bound holds trivially
when $Z<P^{-1}$.

We will need to choose $Z$ as large as possible, given the constraints
in \eqref{23-Z1} and \eqref{23-Z2}.  The choice 
$$
Z=\frac{1}{2}\min\Big\{1, \frac{1}{2c^3fq|z|P^3},
\frac{P}{2q},\max\big\{\frac{q}{c^3fP^3},
q|z|P\big\}\Big\}^{1/3},
$$
is clearly satisfactory. On taking this value in 
\eqref{22-smee} and \eqref{4-mental}, we therefore deduce that
\begin{align*}
|S(\alpha)|^8
& \ll P^{7n+\sigma+1+\ep}
\Big(1+q|z|P^3+qP^{-1}
+q^{-1}\min\big\{P^3,\frac{1}{|z|P}\big\}\Big)^{(n-\sigma-1)/3}\\
& = P^{8n+\ep}
\Big(P^{-3}+q|z|+qP^{-4}
+q^{-1}\min\big\{1,\frac{1}{|z|P^{4}}\big\}\Big)^{(n-\sigma-1)/3},
\end{align*}
whence
\begin{equation}
  \label{eq:chat}
S(\alpha)\ll P^{n+\ep}
\Big(P^{-3}+q|z|\max\{1,\frac{1}{|z|P^{4}}\big\}
+q^{-1}\min\big\{1,\frac{1}{|z|P^{4}}\big\}\Big)^{(n-\sigma-1)/24}.
\end{equation}

We will derive three basic estimates from this bound.  The first
involves the complete exponential sum 
$S_{a,q}$ defined in \eqref{eq:aq}. 
This arises by taking $z=0, P=q$ and $\omega=\chi$ in the 
definition of $S(\alpha)$. An application of \eqref{eq:chat}
immediately gives
\begin{equation}\lab{4-a/q}
S_{a,q}\ll  q^{23n/24+(\sigma+1)/24+\ve},
\end{equation}
for any coprime integers $a,q$ such that $1\leq a \leq q$. Next, we
claim that 
\begin{equation}
  \label{eq:p-2}
S(\al) \ll P^{n+\ep}
(|\alpha|P^{4})^{(\sigma+1-n)/24}, \quad
\mbox{if $|\al|<P^{-2}$}.  
\end{equation}
This is trivial if $|\al|\leq P^{-4}$. If $|\al|>P^{-4}$, then it
follows from \eqref{eq:chat} with $a=0, q=1$ and $\al=z$.
Finally, it is a simple matter to deduce  the following result 
from \eqref{eq:chat}.

\begin{pro}\lab{lem:Sbirch}
Let $a,q, z$ be such that 
$$
1\leq a\leq q, \quad \hcf(a,q)=1, \quad |z|\leq \frac{1}{q^2}. 
$$
Then we have
$$
S(a/q+z)\ll 
P^{n+\ve}
\big(q|z|+q^{-1}|z|^{-1}P^{-4}\big)^{(n-\sigma-1)/24}.
$$
\end{pro}

\section{Estimating $S(\al)$: van der Corput differencing}\label{s:vdc}

We have now come to our own approach for estimating the quartic
exponential sum \eqref{eq:S}, an argument that we have already outlined in \S \ref{section:over}. 
Throughout this section we will assume that
$\omega\in\WW_n$, where $\WW_n$ is defined before the statement of
Lemma~\ref{lem:Iq}. Moreover, we will 
retain the notation \eqref{eq:sigma} for the projective 
dimension of $\sing_\Q(X)$.

Our starting point is the application of van der Corput's method,
which reduces the analysis to a system of cubic exponential sums via
\eqref{vdC1}. As indicated there a comparison of \eqref{vdC1} with
\eqref{weyl1} reveals that the special case $H=P$ of van~der~Corput's method reduces to the first step in
Birch's approach.  Thus we have lost nothing in formulating things
this way, but have gained the considerable advantage that we are now
able to control the size of the parameter $H$ appearing in Propositions \ref{lem:T1}
and \ref{lem:T2}.  For some ranges of $H$ these two propositions,
which use a direct treatment of the cubic sum, are more advantageous
than a result based on Weyl's inequality.

We now proceed to use the results of \S \ref{s:cubic} to 
estimate $T_\h(\alpha)$
directly. We clearly have 
$$
T_\h(\alpha)=\TT_n(\al;F(\x+\h)-F(\x),\omega_{\h},P),
$$
in the notation of  \eqref{eq:Tcubic}.
Moreover, it is not hard to see that
the homogeneous part of $F(\x+\h)-F(\x)$ of maximal degree is
just $\h.\nabla F(\x)$, a cubic form. Furthermore 
$$
\|F(\x+\h)-F(\x)\|_P=\|P^{-3}(F(P\x+\h)-F(P\x))\|\ll H,
$$
since $|\h|\leq H \leq P$.  For $|\h|\leq H$ it is straightforward to check
that $\omega_{\h}\in\WW_n$, where $\omega_{\h}$ is given by
\eqref{wh}.  We now suppose that $\al=a/q+z$ for coprime integers
$a,q$ such that $1\leq a\leq q\leq P^2$, and $z\in\R$ such that 
$|z|\leq q^{-1}P^{-1}$.

Let $\Pi$ denote the set of prime divisors $p$ of $q$, such
that $p^e\| q$, with $e\leq 2$ 
or $e$ odd, and recall the notation
$\Pi_a=\{p\in \Pi: p>a\}$, for any $a\in\N$.
The equation $\h.\nabla F(\x)=0$ defines a
variety $X_{\h}\subseteq \bfP_{\F_v}^{n-1}$, for each 
$v\in\{\infty\}\cup\Pi$.
Let us write
$$
  s_\infty:=s_\infty(X_{\h}), \quad s_p:=s_p(X_\h),
$$
for $p\in\Pi$. In particular, we have $s_p,s_\infty\in 
[-1,n-1]\cap\Z$ and $s_p\geq s_\infty$, for every $p$.
Recall the definition \eqref{eq:qi} of $r_i$,
and write $q=bc^2d$, in the notation of \eqref{eq:bcd}.
Taking $A=1$ in the statement of Proposition \ref{lem:T2}, as
we clearly may, we therefore deduce that
$$
T_\h(\al)\ll
\min_{1+s_\infty\leq \eta \leq n}
q^{-(n-\eta)/2}
\Big(\prod_{i=\eta}^n r_{i}^{(i-\eta)/2}\Big)
P^{n+\ve}W^{n-\eta},
$$
where $W$ is given by \eqref{eq:W}.

We would now like to sum this bound over appropriate values of $\h\in
\Z^n$ in the range $|\h|\leq H$.  Specifically, it follows from
\eqref{vdC1} that there is a factorization $q=r_0\cdots r_n$ such that 
\begin{equation}
  \label{eq:S^2}
|S(\alpha)|^2\ll \frac{P^{2n}}{H^n}\Big(1+ P^\ve
\sum_{s=-1}^{n-1}
\sum_{\h\in\mcal{H}_s}
\min_{1+s\leq \eta\leq n}
\frac{(r_{1+\eta}r_{2+\eta}^2\cdots
  r_n^{n-\eta})^{1/2}}{q^{(n-\eta)/2}} 
W^{n-\eta}\Big),
\end{equation}
where 
$$
\mcal{H}_s:=\Big\{\h\in\Z^n: 
\begin{array}{l}
0<|\h|\leq H,  ~s_\infty(X_\h)=s,\\ 
p\mid r_i
\Rightarrow s_p(X_\h)=i-1,  (0\leq i\leq n)
\end{array}
\Big\}.
$$
Let $m$ be an integer in the range $-1\leq
m\leq n-1$, and let $v\in\{\infty\}\cup\Pi$. It follows from the work of \S
\ref{section:geometry}, and in particular Lemma \ref{lem:dimB}, that 
there exists an affine variety 
$B_{v,m+1}\subseteq \A_{\F_v}^n$ of degree $O(1)$, with
$$
\dim B_{v,m+1} \leq n-m+\sigma_v
$$
for $v\in \{\infty\}\cup \Pi_4$, 
such that $\h\in B_{v,m+1}$ whenever 
$s_{v}(X_{\h})=m$.  Here, $\sigma_v$ is the dimension of the
singular locus of the projective quartic hypersurface $F=0$, viewed
over $\F_v$. On setting $\sigma_\infty=\sigma$, we note that 
there exists an absolute constant $c=O(1)$
such that $\sigma_v=\sigma$ for all $v\in\{\infty\}\cup\Pi_c$.
We may conclude that
$$
\#\mcal{H}_s
\leq \#\big\{\h\in B_{\infty,s+1}\cap \Z^n: |\h|\leq H, ~
[\h]_p \in B_{p,i}~ \forall ~p\mid r_i, ~(0\leq i\leq n)\big\},
$$
where $\dim B_{\infty,s+1}\leq n-s+\sigma$ and $\dim B_{p,i}\leq 
n-i+1+\sigma$, for $0\leq i \leq n$ and $p\in\Pi_c$.

We write, temporarily, 
$$
\rho_i:=\prod_{\colt{p\mid r_i}{p>c}}p
$$ 
for $0\le i\le n$, so
that $\rho_i\gg r_i^{1/2}$ for each index $i$.  Moreover we observe
that
\[i-1=s_p(X_{\ma{h}})\ge s_{\infty}(X_\ma{h})=s\]
for $p\mid r_i$, whence $r_i=1$ for $i\le s$.  We can now apply Lemma
\ref{lem:non-trivial} with $\ell=n-s+\sigma$ and $k_j=n-i+1+\sigma$ for
$p_j\mid r_i$ such that $p_j>c$.  Thus
\begin{align*}
\#\mcal{H}_s
\ll& A(D,n)^{\omega(q)}\Big(H^{n-s+\sigma}\prod_{i=s+1}^n
\rho_i^{s+1-i}\\
& \hspace{3cm}+\omega(q)\sum_{i=s+1}^n H^{n-i+1+\sigma}  
\prod_{j=i+1}^n\rho_j^{i-j}\Big)\\
\ll& q^{\ve}H^{1+\sigma}\sum_{i=s+1}^n
\frac{H^{n-i}}{\prod_{j=i+1}^n\rho_j^{j-i}}\\
\ll& q^{\ve}H^{1+\sigma} \max_{1+s\leq \eta\leq n} 
\frac{H^{n-\eta}}{(r_{1+\eta}r_{2+\eta}^2\cdots r_n^{n-\eta})^{1/2}}.
\end{align*}
It then follows from \eqref{eq:S^2} that
\begin{align*}
|S(\alpha)|^2
&\ll P^{2n+\ve}\sum_{s=-1}^{n-1}
\max_{1+s\leq \eta\leq n}  \frac{H^{1+\sigma-\eta}
}{q^{(n-\eta)/2}} W^{n-\eta}\\
&\ll  P^{2n+\ve}H^{1+\sigma-n} \Big(1+
\frac{H^n}{q^{n/2}} W^{n}\Big).
\end{align*}
On recalling the definition \eqref{eq:W} of
$W$ we therefore deduce the following result.

\begin{pro}\lab{pro:Svdc}
Let $a,q$ be coprime integers  such that 
$1\leq a\leq q\leq P^2$ and $q=bc^2d$, in the notation of \eqref{eq:bcd}.
Let $z\in\R$ such that $|z|\leq q^{-1}P^{-1}$. Then we have
\begin{align*}
S(a/q+z)
\ll \frac{P^{n+\ve}}{H^{(n-1-\sigma)/2}} \Big(1+
\frac{q^{1/2}H}{P} + \sqrt{q|z|H^3P} +
\frac{H}{q^{1/2}}M
\Big)^{n/2},
\end{align*}
where
$$
M:=\min\big\{(c^2d H)^{1/3}, c^{1/2}q^{1/2}(P^{-1/2}+(|z|HP)^{1/4})+ 
c^{5/6}H^{1/6}\big\}.
$$
\end{pro}

\section{Activation of the circle method}\label{s:act}

In this section we recall the apparatus of the Hardy--Littlewood
circle method, as it applies to our problem on quartic hypersurfaces
$X\subset \bfP_\Q^{n-1}$. Let $F\in \Z[x_1,\ldots,x_n]$ be the
underlying quartic form, which we assume to take the shape
\eqref{eq:F}. Our proof of Theorem \ref{main} relies upon using the
circle method to establish an asymptotic
formula for $N_\omega(F;P)$, as 
$P\rightarrow \infty$. Fix once and for all a vector $\x_0\in\R^n$ 
such that $F(\x_0)=0$ and
$\nabla F(\x_0)\neq \ma{0}$. The existence of such a point is
guaranteed since we are assuming in Theorem~\ref{main} that our
hypersurface has a non-singular ad\`elic point, and hence a 
non-singular real point.

We will find it convenient to work with a weight function
that forces us to count points lying very close to $\x_0$. 
For any $\rho\in (0,1]$, we define the function $\omega: \R^n \rightarrow
\R_{\geq 0}$ by 
\begin{equation}
  \label{omega}
\omega(\x):=\gamma(\rho^{-1}|\x-\x_0|),
\end{equation}
where
$$
\gamma(x):=\left\{ 
\begin{array}{ll}
e^{-1/(1-x^2)}, & \mbox{if $|x|<1$},\\
0, & \mbox{if $|x|\geq 1$}.
\end{array}
\right.
$$
Ultimately we will want to work with a value of $\rho$ that is
sufficiently small in some sense, but which still satisfies $1\ll \rho \leq 1$.
It is clear that $\omega$ is infinitely
differentiable, and that it is supported on the region
$|\x-\x_0| \leq \rho$. In particular, we have 
$S(\omega)\leq 1$ provided that $\rho$ and $\x_0$ are sufficiently
small.
Moreover, there exist
constants $c_j>0$ depending only on $j$ and $\rho$ such that
$$
\max\Big\{ \Big|
\frac{\partial^{j_1+\cdots+j_n}\omega (\x)}{\partial^{j_1}x_1\cdots
\partial^{j_n}x_n}\Big|: ~\x \in \R^n, ~j_1+\cdots+j_n=j\Big\}\leq c_j,
$$
for each integer $j \geq 0$. Hence $\omega \in \WW_n$.
Recall the definition \eqref{eq:ss} of the singular series $\ss$, and
define the corresponding singular integral
\begin{equation}
  \label{eq:si}
\mathfrak{I}:=
\int_{-\infty}^{\infty} 
\int_{\R^n} \omega(\x)e(z F(\x))\d\x\d z,
\end{equation}
assuming that it is convergent.
The following result lies at the heart of our proof of Theorem
\ref{main}.

\begin{pro}\label{main'}
Assume that $n-\dim \sing_\Q(X)\geq 42$. Then $\mathfrak{I}$ is
absolutely convergent, and there exists $\delta >0$ such that
$$
N_\omega(F;P)=\mathfrak{S}\mathfrak{I}P^{n-4}+ O\big(P^{n-4-\delta}\big).
$$
\end{pro}

Note that the convergence of $\ss$ is assured by
Theorem \ref{main-ss}. Taking the statement of Proposition
\ref{main'} on faith, we see that in order to complete the proof of Theorem
\ref{main} it will suffice to show that 
$\ss \mathfrak{I}>0$, under the assumption that 
$n-\dim \sing_\Q(X)\geq 42$ and $\Xns(\A_\Q)$ is non-empty. 
The proof that $\ss>0$ follows a standard line of reasoning, as in
\cite[Lemma 7.1]{birch}, and makes use of the fact that $\ss$ is 
absolutely convergent.
Turning to the positivity of $\mathfrak{I}$, we define 
\begin{equation}
  \label{21-si}
\mathfrak{I}(R):=\int_{-R}^{R} 
\int_{\R^n} \omega(\x)e(z F(\x))\d\x\d z.
\end{equation}
for any $R>0$. Then $\mathfrak{I}=\lim_{R\rightarrow \infty}
\mathfrak{I}(R)$. To show that $\mathfrak{I}>0$, it will therefore suffice to
show that $\mathfrak{I}(R)\gg 1$ for sufficiently large values of $R$. 
Performing the integration over $z$, and writing 
$\x=\x_0+\y$, one obtains
\begin{align*}
\mathfrak{I}(R)
&=\int_{\R^n} \omega(\x) \frac{\sin(2\pi R F(\x))}{\pi F(\x)}
\d\x\\
&=\int_{\R^n} \gamma(\rho^{-1}|\y|)  
\frac{\sin(2\pi R F(\x_0+\y))}{\pi F(\x_0+\y)}\d\y.
\end{align*}
The proof that $\mathfrak{I}(R)\gg 1$ is also standard and can be
readily supplied by adapting work of the second author \cite[\S
10] {hb-10} on the corresponding problem for cubic forms.
The only difference lies in the choice of weights used, but this 
does not change the
nature of the proof.  Assume without loss of generality that 
$c_1:=\partial F/\partial x_1(\x_0)\neq 0$.
The need for $\rho>0$ to be sufficiently small 
emerges through an application of the inverse function
theorem. Basically, since $|\y|\leq \rho$, if we write   
$$
z=F(\x_0+\y)=c_1y_1+\cdots+c_ny_n +P_2(\y)+P_3(\y)+P_4(\y)
$$
for forms $P_i$ of degree $i$, then $z\ll \rho$ and we can 
invert this expression to represent $y_1$ as a powers series in 
$z,y_2,\ldots,y_n$, if $\rho$ is sufficiently small. The value of
$\rho$ needed to ensure the validity of such a representation is 
bounded away from zero
in terms of $n$ and $F$ alone.  We refer the reader to \cite{hb-10}
for the remainder of the argument.

As we have already mentioned, Proposition \ref{main'} will
be proved using the circle method. Starting with \eqref{eq:start}, the
idea is to split the interval $[0,1]$ into a
set of major arcs and minor arcs, which are both defined modulo $1$. 
For given $\D>0$,  one takes 
$$
\M_{a,q}(\D):=\Big[\frac{a}{q}- P^{-4+\D}, \frac{a}{q}+ P^{-4+\D}\big]
$$
as major arcs, for $1\le a\le q$ such that $\hcf(a,q)=1$ and 
$q\le P^{\Delta}$.  
It is easily checked that these intervals are disjoint for $\D<4/3$,
which we now assume. Let us write
\begin{equation}
  \label{eq:major}
\M(\D):=\bigcup_{1\leq q\leq P^\D}
\bigcup_{\colt{1\leq a\leq q}{\gcd(a,q)=1}}\M_{a,q}(\D),
\end{equation}
and
$$
\m(\D):=[0,1]\setminus \M(\D)
$$ 
for the corresponding set of minor arcs.  Our treatment of the minor
arc integral $\int_{\m(\D)}S(\alpha)\d\alpha$ will be the focus of
\S\ref{s:minor}, where we will draw together the contents of \S
\ref{s:weyl} and \S \ref{s:vdc}. Next, in \S\ref{s:major} we will
obtain an asymptotic formula for the integral
$\int_{\M(\D)}S(\alpha)\d\alpha$, under suitable hypotheses. These
hypotheses will be validated in \S \ref{s:ss}, during the proof of 
Theorem \ref{main-ss}, which will then complete the proof of
Proposition \ref{main'}.

\section{Treatment of the minor arcs}\label{s:minor}

Recall the definition \eqref{eq:major} of the major arcs $\M(\D)$, for
any $\D<4/3$, and the corresponding set of minor arcs $\m(\D)=[0,1]\setminus
\M(\D)$. Our aim in this section is to establish the following result.

\begin{lem}\label{lem:minor}
Let $n-\dim \sing_\Q(X)\geq  42$.
Then there exists $\del>0$ such that
$$
\int_{\m(\D)}S(\alpha)\d\alpha=O(P^{n-4-\del}),
$$
for any $\D$ in the range $0<\D<4/3$.
\end{lem}

This shows that there is a satisfactory contribution from the minor
arcs in \eqref{eq:start}, when $n-\dim \sing_\Q(X)\geq  42$.
Let $Q\geq 1$. Given any $\al\in[0,1]$, Dirichlet's approximation
theorem allows us to write $\al=a/q+z$, for $a,q\in\N$ and $z\in \R$ such that
$$
1\leq a\leq q, \quad (a,q)=1, \quad q\leq Q, \quad |z|\leq
\frac{1}{qQ}.
$$
In order for such an $\al$ to be contained in the set of minor arcs $\m(\D)$,
it is necessary and sufficient that the inequalities
\begin{equation}
  \label{eq:minor'}
  q\leq P^\D, \quad |z|\leq P^{-4+\D},
\end{equation}
do not both hold.  In our work we will ultimately take the value 
\begin{equation}
  \label{eq:Q}
Q=P^{8/5+\phi},
\end{equation} 
for a small parameter $\phi>0$. This should be 
compared with the value $Q=P^2$ taken by Birch. In effect, using 
van der Corput's method to estimate $S(\al)$
produces a substantially better estimate in the $z$ aspect.
It will be convenient to retain the notation $\sigma=\dim \sing_\Q(X)$,
which was introduced in \eqref{eq:sigma}, and to proceed under 
the assumptions that $n \geq 42+\sigma$ and that $Q$ is contained in the
interval $P^{8/5}\leq Q\leq P^2$.

Given any $R,t\in \R$ such that $0<R\leq Q$ and $0\leq
t\leq(RQ)^{-1}$, we will need to consider 
$$
\Sigma(R,t,\pm):=\sum_{R<q\le 2R}\sum_{\colt{1\leq a\leq q}{\hcf(a,q)=1}}
\int_{t}^{2t}|S(a/q \pm z)|\d z.
$$
Our immediate goal is to obtain conditions on $R$ and $t$ under which
we can establish the existence of $\del>0$ such that
\begin{equation}
  \label{eq:dy_goal}
\Sigma(R,t,\pm)\ll P^{n-4-\del}.  
\end{equation}
If we sum this up over dyadic intervals for $R,t$, it is clear that we
will obtain a satisfactory contribution to the minor arc integral in
Lemma \ref{lem:minor} from the relevant ranges for $R,t$.

Let us begin by considering the overall contribution to $\Sigma(R,t,\pm)$
from those $q$ written in the shape $q=bc^2d$, 
in the notation of \eqref{eq:bcd}, whose factors $b,c,d$ are 
restricted in certain ways. Given $\ma{R}=(R_0,R_1,R_2)\in \R_{\geq 0}$, we let
$\Sigma_{\ma{R}}(R,t,\pm)$ denote the overall contribution to
$\Sigma(R,t,\pm)$ from those $q=bc^2d$ for which
\begin{equation}
  \label{eq:s0s1}
  R_0<b\le 2R_0, \quad   
  R_1<c\le 2R_1, \quad   
  R_2<d\le 2R_2.
\end{equation}
On recalling that $d\mid c$, we note that
$\Sigma_{\ma{R}}(R,t,\pm)=0$ unless
\begin{equation}
  \label{eq:s0s1'}
  R_2\leq 2R_1, \quad R/16< R_0R_1^2R_2\leq 2R, \quad R_i\geq 1/2, 
\end{equation}
for $0\leq i \leq 2$.
The following simple result will be useful in our work.

\begin{lem}\label{lem:count_s0s1}
We have
$$
\sum_{\colt{q=bc^2d}{\mbox{\scriptsize{\eqref{eq:s0s1} holds}}}}
R_0^{i_{0}}R_1^{i_1}R_2^{i_{2}}
\ll R_0^{i_{1}+1}R_1^{i_{1}+1/2}R_2^{i_{2}+1/2},
$$
for any $i_0,i_1,i_2\geq 0$.
\end{lem}

\begin{proof}
We clearly have 
\begin{align*}
\sum_{\colt{q=bc^2d}{\mbox{\scriptsize{\eqref{eq:s0s1} holds}}}}
R_0^{i_{0}}R_1^{i_1}R_2^{i_{2}}
\leq 
\sum_{\colt{b,c,d}{\mbox{\scriptsize{\eqref{eq:bcd}, \eqref{eq:s0s1} hold}}}}
R_0^{i_{0}}R_1^{i_1}R_2^{i_{2}}.
\end{align*}
Recall from \eqref{eq:bcd} that there 
exists a positive integer $d_0$ such that 
$d_0 \mid d$ and $d_0^{-1}d^{-1}c$ is a square-full integer.
Hence, for fixed values of $d$, the number of 
available choices for $b,c$ is
$$
\ll \sum_{d_0 \mid d}R_0 \Big(\frac{R_1}{d_0R_2}\Big)^{1/2}
= R_0R_1^{1/2}R_2^{-1/2}\sum_{d_0 \mid d} \frac{1}{d_0^{1/2}}.
$$ 
On summing over values of $d$, we deduce that the overall number of
choices for $b,c,d$ is
$$
\ll R_0R_1^{1/2}R_2^{-1/2}\sum_{d_0\leq 2R_2}
\frac{1}{d_0^{1/2}}\frac{R_2}{d_0}
\ll R_0R_1^{1/2}R_2^{1/2}.
$$ 
This completes the proof of Lemma \ref{lem:count_s0s1}.
\end{proof}

We are now ready to record the bounds for $\Sigma_{\ma{R}}(R,t,\pm)$ 
that emerge
through our work above.

\begin{lem}\label{lem:min1}
For $t>(RP^2)^{-1}$ we have
$$
\Sigma_{\ma{R}}(R,t,\pm)\ll P^{n+\ve}R_0R_1^{1/2}R_2^{1/2}
(Rt)^{1+(n-\sigma-1)/24},
$$
while for $t\leq(RP^2)^{-1}$ we have
$$
\Sigma_{\ma{R}}(R,t,\pm)\ll
P^{n-(n-\sigma-1)/6+\ve}R_0R_1^{1/2}R_2^{1/2}
(Rt)^{1-(n-\sigma-1)/24}.
$$
\end{lem}
\begin{proof}
This is a straightforward consequence of Proposition \ref{lem:Sbirch}.
Thus we obtain
$$
\Sigma_{\ma{R}}(R,t,\pm)
\ll 
\sum_q\sum_{\colt{1\leq a\leq q}{\hcf(a,q)=1}}
\int_{t}^{2t} P^{n+\ve}
\big(R|z|+R^{-1}|z|^{-1}P^{-4}\big)^{(n-\sigma-1)/24} \d z,
$$
where the summation over  $q$ is over $R<q\le 2R$ such that
\eqref{eq:s0s1} holds. An application of Lemma
\ref{lem:count_s0s1} completes the proof. 
\end{proof}

\begin{lem}\label{lem:min2}
We have
\begin{align*}
\Sigma_{\ma{R}}(R,t,\pm)
\ll&  P^{n/10+9(\sigma+1)/10+\ve}Q^{(n-\sigma-1)/2}\\
&\quad+ RR_0R_1^{1/2}R_2^{1/2}tP^{n+\ve} 
\mu^{(n-\sigma-1)/2},
\end{align*}
where
$$
\mu=\min
\Big\{
\frac{R_1^{5/7}}{R^{3/7}}+R_1^{2/5}t^{1/5}P^{1/5}+\frac{R_1^{1/2}}{P^{1/2}}, ~
\frac{R_1^{1/2}R_2^{1/4}}{R^{3/8}}\Big\}
$$
\end{lem}
\begin{proof}
The proof of Lemma \ref{lem:min2} is based on Proposition \ref{pro:Svdc}. 
An application of this result shows that
for any integer $H$ in the range $1\leq H \leq P$, we have
\begin{align*}
\Sigma_{\ma{R}}(R,t,\pm) 
&\ll
\sum_{q}\frac{R t P^{n+\ve}}{H^{(n-\sigma-1)/2}} 
\Big(1+
\frac{R^{1/2}H}{P} + \sqrt{RtH^3P} +
\frac{HM}{\sqrt{R}}\Big)^{n/2},
\end{align*}
where 
$$
M=\min\{(R_1^2R_2)^{1/3} H^{1/3}, 
R_1^{1/2}R^{1/2}(P^{-1/2}+(tHP)^{1/4})+ R_1^{5/6}H^{1/6}\},
$$
and the summation over  $q$ is over $R<q\le 2R$ such that 
\eqref{eq:s0s1} holds.
We will choose
$$
H=\min\Big\{ \frac{P^{9/5}}{Q}, 
\frac{R^{3/7}}{R_1^{5/7}}, \frac{1}{R_1^{2/5}t^{1/5}P^{1/5}},
\frac{P^{1/2}}{R_1^{1/2}}\Big\}+ \min\Big\{\frac{P^{9/5}}{Q},
\frac{R^{3/8}}{R_1^{1/2}R_2^{1/4}}\Big\}.  
$$
It follows that $H\ll P^{9/5}Q^{-1}$.
Since $\max\{R,P^{8/5}\}\leq Q$ and $t\leq (RQ)^{-1}$ we readily
deduce that $HR^{1/2}P^{-1}\ll 1$ and $RtH^3P\ll 1$.  We then see that
$$
\Sigma_{\ma{R}}(R,t,\pm)
\ll  \sum_{q}
\frac{R t P^{n+\ve}}{H^{(n-\sigma-1)/2}}
\ll
\sum_{q} R t
P^{n+\ve}\Big(\frac{Q}{P^{9/5}}+\mu\Big)^{(n-\sigma-1)/2},
$$
for this choice of $H$, where $\mu$ is as in the statement of the lemma.
The contribution from the term involving $Q/P^{9/5}$ is 
$$
\ll P^{n/10+9(\sigma+1)/10+\ve}Q^{(n-\sigma-1)/2}.
$$
We complete the proof of Lemma~\ref{lem:min2} via an 
application of Lemma~\ref{lem:count_s0s1} in the above.
\end{proof}

We clearly have
$$
\Sigma(R,t,\pm)\ll P^\ve \max_{R_0,R_1,R_2} \Sigma_{\ma{R}}(R,t,\pm),
$$
where the maximum is over all values of $\ma{R}=(R_0,R_1,R_2)\in
\R_{\geq 0}^3$ such that \eqref{eq:s0s1'} holds.
Let $\D>0$ and recall the inequalities \eqref{eq:minor'}. 
We want to show that \eqref{eq:dy_goal} holds unless
\begin{equation}
  \label{eq:dy_minor}
\mbox{$2R\leq P^\D$ ~and~ $ 2t\leq P^{-4+\D}$}.  
\end{equation}
We proceed by considering the two basic ranges for $t$ that will
emerge through our application of Lemma \ref{lem:min1}.

Suppose first that $t$ lies in the range
\begin{equation}
  \label{eq:trange}
(RP^2)^{-1}<t\leq (RQ)^{-1},
\end{equation}
and write $M_0=P^{n/10+9(\sigma+1)/10+\ve}Q^{(n-\sigma-1)/2}.$
Then Lemmas \ref{lem:min1} and \ref{lem:min2} give
\begin{align}\label{eq:t1}
\Sigma(R,t,\pm)
&\ll M_0+
\max_{\ma{R}}RR_0R_1^{1/2}R_2^{1/2}tP^{n+\ve}
\min\Big\{R^{1/24}t^{1/24}\,,\,\mu^{1/2}\Big\}^{n-\sigma-1}
\nonumber\\ 
&\ll  M_0+
\max_{\ma{R}}\frac{R_0R_1^{1/2}R_2^{1/2}P^{n+\ve}}{Q}
(M_1+M_2+M_3)^{n-\sigma-1},
\end{align}
where we have set
\begin{align*}
M_1&=\frac{1}{2}\min\Big\{\frac{1}{Q^{1/24}}, 
\frac{R_1^{1/4}R_2^{1/8}}{R^{3/16}}, 
\frac{R_1^{5/14}}{R^{3/14}}\Big\},\\
M_2&=\frac{1}{2}\min\Big\{\frac{1}{Q^{1/24}}, 
\frac{R_1^{1/4}R_2^{1/8}}{R^{3/16}}, 
\frac{R_1^{1/5}P^{1/10}}{R^{1/10}Q^{1/10}}
\Big\},\\
M_3&=\frac{1}{2}\min\Big\{\frac{1}{Q^{1/24}}, 
\frac{R_1^{1/4}R_2^{1/8}}{R^{3/16}}, 
\frac{R_1^{1/4}}{P^{1/4}}\Big\}.
\end{align*}
Let us write $\Sigma_i$ for the overall contribution to
$\Sigma(R,t,\pm)$ from the term involving $M_i$, for $i=1,2,3$. Then
we have
\begin{equation}
  \label{eq:123}
\Sigma(R,t,\pm)
\ll M_0+ \Sigma_1+\Sigma_2+\Sigma_3.
\end{equation}
We now wish to show that each of the terms $M_0, \Sigma_1, 
\Sigma_2, \Sigma_3$ is
$o(P^{n-4})$, when $Q$ takes the value \eqref{eq:Q} for a suitable
choice of $\phi\in(0,2/5)$.  
Beginning with $M_0$, our choice of $Q$ yields
$$
M_0\leq P^{9n/10+(\sigma+1)/10+\ve+\phi n/2}.
$$
The exponent here is strictly less than $n-4$ for $n\geq 42+\sigma$,
provided that $\ve$ is sufficiently small and
$\phi<n^{-1}(1/5-2\ve)$.

Turning to the terms $\Sigma_1,\Sigma_2,\Sigma_3$, we note that $M_i\leq 1$ for
$i=1,2,3$, since \eqref{eq:s0s1'} implies that 
\begin{equation}
  \label{eq:expand}
R_1^{1/4}R_2^{1/8}R^{-3/16}\leq 2. 
\end{equation}
Hence $M_i^{n-\sigma-1}$ is a decreasing function of $n$. Thus it
will suffice to show that $\Sigma_i=o(P^{n-4})$ at the value 
$n=42+\sigma$, which we
now assume.  In what follows, we will make frequent use of the
inequality 
$$
\min\{A,B,C\}\leq A^\al B^\be C^\gamma,
$$
for any $\al,\be,\gamma\geq 0$ such that $\al+\be+\gamma=1$.
Let us begin by considering the term involving $M_1$. Note that
\begin{align*}
M_1^{n-\sigma-1}&\leq  
\Big(\frac{1}{Q^{1/24}}\Big)^{36} 
\Big(\frac{R_1^{1/4}R_2^{1/8}}{R^{3/16}}\Big)^{8/3}
\Big(\frac{R_1^{5/14}}{R^{3/14}}\Big)^{7/3}=
\frac{R_1^{3/2}R_2^{1/3}}{Q^{3/2}R},
\end{align*}
whence
\begin{align*}
\Sigma_1&\ll 
\max_{\ma{R}}\frac{R_0R_1^{2}R_2^{5/6}P^{n+\ve}}{Q^{5/2}R}
\ll \frac{P^{n+\ve}}{Q^{5/2}}=P^{n-4-5\phi/2+\ve}.
\end{align*}
This is clearly satisfactory provided that 
$\phi>2\ve/5$. Next, we observe that 
$$
M_2^{n-\sigma-1}\leq  
\Big(\frac{1}{Q^{1/24}}\Big)^{27} 
\Big(\frac{R_1^{1/4}R_2^{1/8}}{R^{3/16}}\Big)^{4}
\Big(\frac{R_1^{1/5}P^{1/10}}{R^{1/10}Q^{1/10}}\Big)^{10}=
\frac{R_1^{3}R_2^{1/2}P}{Q^{17/8}R^{7/4}}.
$$
Thus
$$
\Sigma_2\ll 
\max_{\ma{R}}\frac{R_0 R_1^{7/2}R_2 P^{n+1+\ve}}{Q^{25/8}R^{7/4}}
\ll \frac{P^{n+1+\ve}}{Q^{25/8}}=P^{n-4-25\phi/8+\ve}. 
$$
This too is satisfactory provided that $\phi>8\ve/25$.
Finally, we note that 
$$
M_3^{n-\sigma-1}\leq  
\Big(\frac{1}{Q^{1/24}}\Big)^{31} 
\Big(\frac{R_1^{1/4}R_2^{1/8}}{R^{3/16}}\Big)^{8}
\Big(\frac{R_1^{1/4}}{P^{1/4}}\Big)^{2}=
\frac{R_1^{5/2}R_2}{P^{1/2}Q^{31/24}R^{3/2}},
$$
whence
$$
\Sigma_3\ll 
\max_{\ma{R}}\frac{R_0 R_1^{3}R_2^{3/2} P^{n-1/2+\ve}}{Q^{55/24}R^{3/2}}
\ll \frac{P^{n-1/2+\ve}}{Q^{55/24}}=P^{n-25/6-55\phi/24+\ve}.
$$
This is satisfactory for any $\phi\geq 0$, provided that $\ve$ is
sufficiently small.  Feeding these results into 
\eqref{eq:123}, we therefore see
that \eqref{eq:dy_goal} holds for any $t$ in the range
\eqref{eq:trange}, 
with the choice \eqref{eq:Q} for $Q$,
provided that
$$
\frac{2\ve}{5}<\phi <\frac{1}{5n}-\frac{2\ve}{n}.
$$
Such an interval is clearly non-empty for $\ve$ sufficiently small.

In order to complete the treatment of the minor arcs it remains to
produce a similar bound for $t$ in the range
\begin{equation}
  \label{eq:trange'}
t\leq (RP^2)^{-1},
\end{equation}
still under the assumption that \eqref{eq:dy_minor} does not hold.
In this case we get a bound similar to \eqref{eq:t1}, but with 
the term $P^{n+\ve}/Q$ replaced by $RtP^{n+\ve}$, and with
$Q^{-1/24}$ replaced by $(RtP^4)^{-1/24}$ in the definitions of
$M_1,M_2,M_3.$
In particular \eqref{eq:expand} still ensures that $M_i \leq 1$ for $i=1,2,3.$
In our estimation of $\Sigma_i$, for $i=1,2,3$, we made use of the
observation that 
$M_i^{n-\sigma-1}\leq A^\al B^\be C^\gamma$ for 
$n= 42+\sigma$, where $A,B,C$ are the
three terms in the definition of $M_i$ and $\al+\be+\gamma=41$. 
Thus, if $\al_i=36$, 27, or 31 denotes the exponent of $A$ that was chosen to
estimate $\Sigma_i$, we see that in order to estimate the contribution from
the new range for $t$, it will suffice to multiply
the final answer by
$$
E_{\al_i}:=\Big(\frac{Q}{Rt P^4}\Big)^{\al_i/24} RtQ
$$
in each case.  Thus for $t$ in the range
\eqref{eq:trange'} and $Q$ given by \eqref{eq:Q}, for a suitable value
of $\phi\in[0,2/5]$, it suffices the check that $E_{\al}=O
(1)$ at $\al=27$ and $36$, unless \eqref{eq:dy_minor} holds.
This we proceed to do.

Note first that 
\begin{align*}
E_{27} =\frac{Q^{17/8}}{R^{1/8}t^{1/8}P^{9/2}},\quad
E_{36} =\frac{Q^{5/2}}{R^{1/2}t^{1/2}P^6}.
\end{align*}
Hence these terms are satisfactory if
$$
t\gg \max\Big\{\frac{Q^{17}}{RP^{36}}\,,\,\frac{Q^5}{RP^{12}}\Big\}=
\frac{Q^5}{RP^{12}}=\frac{1}{RP^{4-5\phi}}. 
$$
It remains to deal with the possibility that 
$t\ll R^{-1}P^{-4+5\phi}$. Let us suppose that $\phi<\D/5$, so that in
particular we have $2t\leq P^{-4+\D}$ for $P$ sufficiently large.
To deal with $t\ll R^{-1}P^{-4+5\phi}$, we may therefore assume 
that $2R> {P^\D}$, since we are supposing that \eqref{eq:dy_minor}
does not hold. An application of Lemma
\ref{lem:min2} now gives
\begin{align*}
\Sigma(R,t,\pm)
&\ll M_0+ \max_{\ma{R}}
RR_0R_1^{1/2}R_2^{1/2}tP^{n+\ve} 
\Big(\frac{R_1^{1/2}R_2^{1/4}}{R^{3/8}}\Big)^{(n-\sigma-1)/2}\\
&\ll M_0+ \max_{\ma{R}}
R_0R_1^{1/2}R_2^{1/2}P^{n-4+5\phi+\ve} 
\Big(\frac{1}{R_0^{3}R_1^{2}R_2}\Big)^{41/16}\\
&= M_0+ \max_{\ma{R}}
R_0^{-107/16}R_1^{-37/8}R_2^{-33/16}P^{n-4+5\phi+\ve} \\
&\ll M_0+ R^{-33/16} P^{n-4+5\phi+\ve} \\
&\ll M_0+ P^{n-4-33\D/16+5\phi+\ve}
\end{align*}
for $n\geq 42+\sigma$, 
where the maximum is over all values of $\ma{R}\in
\R_{\geq 0}^3$ such that \eqref{eq:s0s1'} holds.
We previously gave a satisfactory treatment of $M_0$ for $2\ve/5<\phi< 
n^{-1}(1/5-2\ve)$. Thus \eqref{eq:dy_goal} holds for any $t$ in the range
\eqref{eq:trange'}, 
provided that \eqref{eq:dy_minor} does not hold and assuming that
$\phi$ lies in the range 
$$
\frac{2\ve}{5}<\phi <\min\Big\{\frac{1}{5n}-\frac{2\ve}{n},\frac{\D}{5}\Big\}.
$$
Such an interval is clearly non-empty for any $\D>0$, if 
$\ve$ is chosen to be sufficiently small.
This completes the proof of Lemma \ref{lem:minor}.

\section{Treatment of the major arcs}\label{s:major}

Recall the definition \eqref{eq:aq} of the complete exponential sums
$S_{a,q}$, for any coprime integers $a,q$ such that $1\leq a\leq q$, and set
\begin{equation}\label{21-sing}
\ss(R):=\sum_{q\leq R} \frac{1}{q^n}\sum_{\colt{a=1}{\hcf(a,q)=1}}^q
S_{a,q},
\end{equation}
for any $R>1$. Then $\ss=\lim_{R\rightarrow \infty}\ss(R)$
in \eqref{eq:ss}. Recall the definition \eqref{eq:si} of $\mathfrak{I}$.
The aim of this section is to establish the
following result.

\begin{lem}\label{lem:major}
Let $n-\dim \sing_\Q(X)\geq  26$. Suppose that 
$\ss$ is absolutely convergent, and satisfies the
estimate
\begin{equation}\lab{27-record}
\mathfrak{S}(R)=\mathfrak{S}+O_\phi(R^{-\phi}),
\end{equation}
for some $\phi>0$.  Then $\mathfrak{I}$ is absolutely 
convergent, and there exists $\del>0$
such that
$$
\int_{\M(\D)}S(\alpha)\d\alpha=\ss \mathfrak{I}P^{n-4}
+O_{\phi}(P^{n-4-\del}),
$$
for any $\D$ in the range $0<\Delta<1/5$. 
\end{lem}

The statement of Proposition \ref{main'} now follows from 
Lemmas \ref{lem:minor} and \ref{lem:major},  together 
with Theorem \ref{main-ss}, under
the further assumption that \eqref{27-record} holds for some
$\phi>0$. This latter estimate will be established in \S \ref{s:ss},
for $n$ in the range $n-\dim\sing_\Q (X)\geq 42$.

Turning to the proof of Lemma \ref{lem:major}, let 
$\al\in \M_{a,q}(\D)$, with $\al=a/q+z$. Furthermore, let $\omega$ be
given by \eqref{omega}. Following the lines of 
\cite[Lemma~5.1]{birch}, we will show that
\begin{equation}
  \label{eq:train1}
S(\al)=q^{-n}P^n S_{a,q} I(z P^4) +O(P^{n-1+2\D}),  
\end{equation}
where $S_{a,q}$ is given  by \eqref{eq:aq} and
$$
I(\gamma):=\int_{\R^n} \omega(\x)e(\gamma F(\x)) \d\x,
$$
for any $\gamma \in \R$.  To see this we write $\x=\y+q\z$, where $\y$
runs over a complete set of residues modulo $q$, giving
\begin{equation}
  \label{eq:train2}
  S(\al)=\sum_{\y \bmod{q}}e_q\big(aF(\y)\big)
\sum_{\z\in\Z^n}f(\z),
\end{equation}
where 
$$
f(\z)=\omega\Big(\frac{\y+q\z}{P}\Big)e\big(z F(\y+q\z)\big).
$$
Note that $q\leq P^\D$ by assumption. 
We now want to replace the discrete variable $\z$ by a continuous
one, and to replace the summation over $\z$ by an integration. 
For this purpose it will suffice to use a rather crude estimate.
If $\x\in [0,1]^n$ then $f(\z+\x)=f(\z)+O(\max_{\ma{u}\in [0,1]^n} 
|\nabla
f(\z+\ma{u})|)$. Hence
\begin{align*}
\Big|\int_{\R^n}f(\z)\d\z-\sum_{\z\in\Z^n}f(\z)\Big|
&\ll \meas (\mcal{S})
\max_{\z\in \mcal{S}}|\nabla f(\z)|\\
&\ll \Big(\frac{P}{q}\Big)^n \big(q/P+q|z|P^3)\\
&\ll |z|q^{1-n}P^{n+3}+q^{1-n}P^{n-1},
\end{align*}
where $\mcal{S}$ is an $n$-dimensional cube with sides of order
$1+P/q\leq 2P/q$. Hence
\begin{align*}
\sum_{\z\in\Z^n}f(\z)&=\frac{P^n}{q^{n}} \int_{\R^n} 
\omega(\x)e\big(z P^4F(\x)\big)\d \x
+O\big(|z|q^{1-n}P^{n+3}+q^{1-n}P^{n-1} \big),
\end{align*}
on making the change of variables $P\x=\y+q\z$. Substituting this
into \eqref{eq:train2}, we therefore deduce that
\begin{equation}
  \label{eq:train3}
  S(\al)=q^{-n}P^n S_{a,q} I(z P^4) +O(|z|qP^{n+3}+qP^{n-1}).
\end{equation}
This completes the proof of \eqref{eq:train1}, since
$|z|\leq P^{-4+\D}$ and $q\leq P^\D$ on the major arcs.

Using \eqref{eq:train1}, and noting that the major arcs have measure
$O(P^{-4+3\D})$, it is now a trivial matter to deduce that
\begin{equation}\lab{21-maj}
\int_{\M(\D)} S(\al)\d\al=
P^{n-4}\mathfrak{S}(P^\D)\mathfrak{I}(P^\D) +O(
P^{n-5+5\D}),
\end{equation}
where $\mathfrak{S}(P^\D)$ is given by \eqref{21-sing}, and
$\mathfrak{I}(P^\D)$ is given by \eqref{21-si}.
Recall the definition \eqref{eq:si} of the singular 
integral $\mathfrak{I}$. Then
we clearly have $\mathfrak{I}=\lim_{R\rightarrow \infty}
\mathfrak{I}(R)$, if this
limit exists.  Write $\sigma=\dim\sing_\Q(X)$, as in 
\eqref{eq:sigma}. 
We now estimate $I(\gamma)$.

\begin{lem}\lab{21-int1}
We have $I(\gamma)\ll \min\{1,|\gamma|^{(\sigma+1-n)/24+\ve} \}$.
\end{lem}

\begin{proof}
The estimate $I(\gamma)\ll 1$ is trivial. In proving the second
estimate we may clearly assume that $|\gamma|>1$. 
Taking $a=0$ and $q=1$ in \eqref{eq:train3}, we deduce that
$$
S(\al)=P^n I(\al P^4) +O\big((|\al| P^4+1)P^{n-1}\big),
$$
for any $P\geq 1$.  On the other hand, assuming that $|\al|<P^{-2}$,
\eqref{eq:p-2} gives
$$
S(\alpha)  \ll 
P^{n+\ve}(|\al|P^4)^{(\sigma+1-n)/24}. 
$$
Writing $\al P^4=\gamma$, we may combine these estimates to obtain
$$
I(\gamma)\ll 
|\gamma|^{(\sigma+1-n)/24}P^{\ve}
+|\gamma| P^{-1},
$$
when $|\gamma|<P^{2}$.  Finally we observe that $I(\gamma)$ is
independent of $P$.  Thus we are free to choose
$P=|\gamma|^{(24+n-\sigma-1)/24}$, from which the second estimate 
of Lemma \ref{21-int1} follows.
\end{proof}

Suppose that $n\geq 26+\sigma$.
It now follows from Lemma \ref{21-int1} that
\begin{align*}
\mathfrak{I}-\mathfrak{I}(R)
=\int_{|\gamma|\geq R} I(\gamma) \d\gamma
&\ll
\int_{R}^{\infty} \min\{1,\gamma^{(\sigma+1-n)/24+\ve}\} \d\gamma\\
&\ll R^{(25+\sigma-n)/24+\ve}.
\end{align*}
This shows in particular that $\mathfrak{I}$ is absolutely convergent
for $n\geq 26+\sigma$, as claimed
in the statement of Lemma \ref{lem:major}. 
In particular, $\mathfrak{I}(P^\D)\ll 1$. Hence, under the assumption
that \eqref{27-record} holds and $n\geq 26+\sigma$, it follows
from \eqref{21-maj} that 
\begin{align*}
\int_{\M(\D)} S(\al)\d\al
&=
\ss P^{n-4}\mathfrak{I}(P^\D)
+O_\phi\big(P^{n-5+5\D}+P^{n-4-\D\phi}\big)\\
&=
\ss \mathfrak{I}P^{n-4}
+O_\phi\big(P^{n-5+5\D}+P^{n-4-\D\phi}+P^{n-4-\D/24+\ve}\big),
\end{align*}
with $\phi>0$.
We therefore obtain the statement of Lemma \ref{lem:major} by choosing
$\D$ such that $0<\D<1/5$, and taking $\ve>0$ to be
sufficiently small.

\section{The singular series}\label{s:ss}

Recall the definition \eqref{eq:ss} of the singular series $\ss$, and
that of the complete exponential sum \eqref{eq:aq}.
It is well-known and easy to check that the summands 
$$
q^{-n}\sum_{\colt{1\leq a \leq q}{\hcf(a,q)=1}} S_{a,q}
$$
are multiplicative functions of $q$. 
It follows that the singular series $\ss$ is absolutely 
convergent if and only if the
product
$
\prod_p (1+\sum_{k=1}^\infty a_p(k))
$
is, where 
$$
a_{p}(k):=
p^{-kn}\sum_{\colt{1\leq a\leq p^k}{\gcd(a,p)=1}} |S_{a,p^k}|.
$$
This product is absolutely convergent if and only if the sum
$
\sum_p \sum_{k=1}^\infty a_{p}(k)
$ 
is convergent.
Now it follows from \eqref{4-a/q} that
\begin{equation}
  \label{eq:ap-1}
  a_p(k)\ll p^{k(1+(\sigma+1)/24-n/24)+\ve},
\end{equation}
for any $k\geq 1$, where $\sigma$ is given by \eqref{eq:sigma} as usual. 
This establishes the 
absolute convergence of $\ss$ for $n\geq 50+\sigma$, which is enough
for Birch's result.

Let us assume henceforth that $n\geq 27+\sigma$. Then \eqref{eq:ap-1}
yields 
$$
\sum_p\sum_{k\geq 24} a_p(k)\ll \sum_p p^{25+\sigma-n+\ve} \leq
\sum_{m=1}^\infty
\frac{1}{m^{2-\ve}}\ll 1
$$
if $\ve$ is sufficiently small. To handle the contribution 
from the $a_p(k)$ for
$2\leq k<24$, we will employ the following simple result

\begin{lem}\label{k>1}
We have 
$$
S_{a,p^k} \ll_k p^{(k-1)n+\sigma+1}
$$
for $k\geq 2$.
\end{lem}

\begin{proof} 
Let $k\geq 2$ and write $\x=\y+p^{k-1} \z$. Then it follows that
\begin{align*}
S_{a,p^k}&= \sum_{\y\bmod{p^{k-1}}} e_{p^k}(aF(\y)) \sum_{\z \bmod{p}}
e_p(a\z.\nabla F(\y) ),
\end{align*}
whence
\begin{align*}
|S_{a,p^k}| 
&\leq
p^n \#\{\y \bmod{p^{k-1}}: p\mid \nabla F(\y)\}\\
&=
p^{(k-1)n} \#\{\y \bmod{p}: p\mid \nabla F(\y)\}\\
&\ll
p^{(k-1)n+\sigma+1} 
\end{align*}
for $p\gg 1$. Since $S_{a,p^k}\ll_k 1$ for the remaining 
values of $p$, this therefore suffices for the proof of the lemma.
\end{proof}

It follows from Lemma \ref{k>1} that
$$
\sum_p \sum_{k=2}^{23} a_p(k) \leq 
\sum_p \sum_{k=2}^{23} p^{k(1-n)}.p^{(k-1)n+\sigma+1} \ll 
\sum_p p^{24+\sigma-n}\ll 1,
$$
still under the assumption that $n\geq 27+\sigma$.  It remains to deal
with the sum $\sum_p a_p(1)$. For this we can apply Lemma
\ref{deligne} to conclude that
$$
\sum_p a_p(1)\ll \sum_{p}p^{1-n/2+(\sigma+1)/2}\ll 1,
$$
which therefore completes the proof of Theorem \ref{main-ss}.

It remains to establish the estimate in \eqref{27-record} for a
suitable $\phi>0$, where $\ss(R)$ is given by \eqref{21-sing}.
This we will do under the assumption that $n\geq 42+\sigma$, as we
clearly may.  Recall from Lemma \ref{lem:T=mult} that 
$S_{a,q}$ is multiplicative in $q$. 
Let us write $q=uv$, where $u$ is the square-free part of $q$. 
Then it follows from Lemma \ref{deligne} that
$
S_{a,u}\ll u^{(n+\sigma+1)/2+\ve}.
$
Once combined with \eqref{4-a/q} we deduce that
\begin{align*}
\big|\ss-\ss(R)\big| 
&\ll \sum_{q=uv>R} u^{1-n/2+(\sigma+1)/2+\ve}
v^{1-n/24+(\sigma+1)/24+\ve}\\
&\ll \sum_{q=uv>R} u^{-39/2+\ve} v^{-17/24+\ve}\\
&\ll R^{-5/24+2\ve}\sum_{q=uv>R} u^{-2} v^{-1/2-\ve}\\
&\ll R^{-5/24+2\ve}\sum_{u,v=1}^\infty u^{-2} v^{-1/2-\ve},
\end{align*}
since $n\geq 42+\sigma$.
Now the number of square-full integers $v\in (V,2V]$ is $O(V^{1/2})$, whence 
the sum over $v$ is convergent, as, of course, is the sum over $u$.
Thus $\big|\ss-\ss(R)\big| \ll R^{-5/24+2\ve}$. 
This therefore completes the proof of
\eqref{27-record}, in which one can take any $\phi\in(0,5/24)$.

\end{document}